\documentclass[a4paper, 12pt]{article}
\usepackage{amssymb}
\usepackage{stmaryrd}
\usepackage{amsmath}
\usepackage{amsthm}
\usepackage{tipa}
\usepackage{faktor}
\usepackage{dsfont}
\usepackage{mathtools}
\usepackage{lipsum}
\usepackage{soul}
\usepackage{tikz}
\usepackage{amscd}
\usepackage{t1enc}
\usepackage[mathscr]{eucal}
\usepackage{indentfirst}
\usepackage{epic}
\usepackage{epstopdf}
\usepackage{textcomp}
\usepackage{cleveref}
\usepackage{setspace}
\usepackage{float}
\usepackage[bottom]{footmisc}
\usepackage{authblk}
\usepackage{fancyhdr}
\renewcommand{\qed}{$\blacksquare$}
\usepackage[margin=2.9cm]{geometry}

\theoremstyle{plain}
\newtheorem{theorem}{Theorem}
\newtheorem*{theorem*}{Theorem}
\newtheorem{lemma}{Lemma}

\newtheorem{corollary}{Corollary}
\newtheorem{remark}{Remark}

\newtheorem*{conjecture*}{Conjecture}
\usepackage{graphicx}
\graphicspath{{./symmetryimages/}}
\usepackage{caption}

\newcommand\norm[1]{\left\lVert#1\right\rVert}

\newcommand\blfootnote[1]{%
	\begingroup
	\renewcommand\thefootnote{}\footnote{#1}%
	\addtocounter{footnote}{-1}%
	\endgroup
}

\newtheorem*{innerthm}{\innerthmname}
\newcommand{\innerthmname}{}
\newenvironment{statement}[1]
{\renewcommand{\innerthmname}{#1}\innerthm}
{\endinnerthm}

\numberwithin{equation}{section}

\begin{document}
	
	\title{\textbf{On Proinov's lower bound for the diaphony}}
	
	\author{\Large{Nathan Kirk} \\ \small{Queen's University Belfast} \\ \textit{Email: nkirk09@qub.ac.uk}}

	\date{}
\maketitle

\begin{abstract}
\noindent
In 1986, Proinov published an explicit lower bound for the diaphony of 
both finite and infinite sequences of points contained in the $d-$dimensional unit cube \cite{Proinov 1}. However, his widely cited paper does not contain the proof of this result but simply states that this will appear elsewhere. To the best of our knowledge, this proof was so far only available in a monograph of Proinov written in Bulgarian \cite{Proinov 2}. The first contribution of our paper is to give a self contained version of Proinov's proof in English. Along the way, we improve the explicit asymptotic constants implementing recent, and corrected results of Hinrichs \& Markhasin \cite{Aicke}, and Hinrichs \& Larcher \cite{HinLar}. [The corrections are due to a note in \cite{HinLar}]. Finally, as a main result, we use the method of Proinov to derive an explicit lower bound for the dyadic diaphony of finite and infinite sequences in a similar fashion.    


\end{abstract}

\section{Introduction}

The beginnings of the theory of uniform distribution modulo one can be attributed to the work of H. Weyl \cite{WeylCrit} of 1916. J. Van der Corput \cite{vandercor1, vandercor2} later conjectured that no sequence can be, in some sense, too evenly distributed. In 1954, K. Roth \cite{Roth1954} improved on the thoughts of Van der Corput publishing a celebrated sharp lower bound for the $\mathcal{L}_{2}-$discrepancy of an $d-$dimensional finite sequence, $\sigma_N$. In particular,
\begin{equation*}
\mathcal{L}_{2, N}(\sigma_N) \geq c(d) \frac{(\log N)^{\frac{d-1}{2}}}{N}
\end{equation*}
where $c(d)$ is a constant dependent only upon the dimension. The result of Roth has specific importance throughout this paper. For a more detailed and comprehensive history of the beginnings and development of the quantitive measures of uniform distribution theory, we refer the reader to the survey \cite{Bilyk}.

\vspace{1mm}
Motivated by the heavy influence of trigonometric summations in Weyl's Criterion for uniform distribution modulo one and the inequality of Erd\"{o}s-Tur\'{a}n \cite{ErdosTuran, KuipNied}, P. Zinterhof proposed a new measure of irregularity of distribution in \cite{Zinterhof} which he named, \textit{diaphony}, denoted throughout by $F_N$. Similar to the above result for the $\mathcal{L}_2-$discrepancy, in 1986 P. Proinov published results \cite{Proinov 1} allowing one to calculate exact lower bounds for the diaphony of arbitrary $d-$dimensional sequences.

\vspace{1mm}
Of particular concern to the author is a simple corollary of Proinov's work concerning the lower bound of one-dimensional sequences contained in the unit interval. It is known \cite{Proinov 1} and will be shown in this paper, that for an infinite one-dimensional sequence $\sigma$,
\begin{equation*} \label{diaphonylb}
F_N(\sigma) \geq c \cdot \frac{\sqrt{\log N}}{N} \tag{$*$}
\end{equation*}
holds for infinitely many $N$, where $c > 0$ is an absolute constant. It is therefore natural to consider, what is the largest value of $c$ for which (\ref{diaphonylb}) holds for all one-dimensional sequences $\sigma$ for infinitely many $N$? To investigate, we define the asymptotic constant for the diaphony of an infinite one-dimensional sequence $\sigma$,
\begin{equation*}
f(\sigma) := \limsup_{N \rightarrow \infty} \frac{N F_N(\sigma)}{\sqrt{\log N}}
\end{equation*}
and denote by,
\begin{equation*}
f^* := \inf_{\sigma} f(\sigma)
\end{equation*}
the \textit{one-dimensional diaphony constant}. That is, $f^*$ is the supremum over all $c$ such that (\ref{diaphonylb}) holds. Study in areas of the same flavour have appeared recently in the form of asympototic constants of the corresponding notions of (star and extreme) discrepancy \cite{Larch1, Larch2, LarchPuch, Pausinger}. Returning to our motivation, in his 1986 paper Proinov states a lower bound for $f^*$. This paper is widely cited however the proofs of several of the results are not included in the text and instead, they are simply said to appear elsewhere. Therefore, the first aim of this paper is to make these proofs accessible. Further to collating these hidden proofs, and due to recent results improving lower bounds of the $\mathcal{L}_2-$discrepancy in \cite{Aicke} and \cite{HinLar}, we update and improve the results concerning lower bounds of the diaphony of general $d-$dimensional sequences and state the up-to-date one-dimensional diaphony constant. 

\vspace{1mm}
As discussed above, the concept of the diaphony is based on the trigonometric function system. However, introduced by Hellekalek and Leeb in \cite{Hellekalek}, another notion of diaphony exists based on the (dyadic) Walsh function system.\footnote{J. Walsh published his namesake function system in 1923, \cite{Walsh}. } This is aptly named, \textit{dyadic diaphony} and denoted throughout by $F_{2, N}$.\footnote{It was found that there exists an innate relationship between the function system that is chosen and the type of constructions of sequences that can be analysed with the corresponding Weyl summations. For example, the trigonometric function system is well suited to study lattice point sequences and in this instance, the Walsh function system is better suited to analyse digital nets and sequences, \cite{Larch3}. } It is already known \cite{Pill1} that for the dyadic diaphony, 
\begin{equation*}
F_{2, N}(\sigma_N) \geq \bar{c}(d) \frac{(\log N)^{\frac{d-1}{2}}}{N}
\end{equation*}
where $\sigma_N$ is a finite sequence contained in the $d-$dimensional unit cube, and $\bar{c}(d)$ is a constant dependent only upon the dimension. In this paper, after understanding Proinov's methods in the case of the classical diaphony, we move in the latter stages to use these same techniques in the setting of the dyadic diaphony. In doing so, we arrive at analogous explicit lower bounds for the dyadic diaphony and hence finish by stating an equivalent lower bound for the \textit{one-dimensional dyadic diaphony constant}, 
\begin{equation*}
f^*_2 :=\inf_{\sigma} \limsup_{N \rightarrow \infty} \frac{N F_{2, N}(\sigma)}{\sqrt{\log N}}.
\end{equation*}

In what follows, Section 2.1 gives the necessary preliminaries which allow the statement of Proinov's Theorems in Section 2.2. We proceed to give the means in which we can state the updated constant for the diaphony and a new constant for the dyadic diaphony in Sections 2.3 and 2.4 respectively. Section 3.1 contains a high level overview of the proof of Proinov, while Section 3.2 follows to give full, detailed proofs. Lastly, Section 4 gives a proof for the main result in the derivation of the explicit lower bound for the dyadic diaphony.

\section{Statement of Results}

\subsection{Preliminaries and Notation}

\paragraph{Discrepancy.}
In this paper we are concerned with the distribution of points in the $d-$dimensional unit cube, $[0, 1)^d$. Let $\sigma_N = (\mathbf{a}_i)_{i=1}^N$ be a finite sequence of points contained in $[0, 1)^d$. For any point $\boldsymbol{\gamma} = (\gamma_1, \gamma_2, ..., \gamma_d) \in [0, 1)^d$ define the discrepancy function as,
\begin{equation*}
g \big( [\mathbf{0}, \boldsymbol{\gamma}), \sigma_N, N \big) \coloneqq \frac{1}{N} \sum_{i=1}^{N} \chi_{\boldsymbol{\gamma}}(\mathbf{a}_i) - \lambda_d \big( [\mathbf{0}, \boldsymbol{\gamma}) \big)
\end{equation*}
where $\chi_{\boldsymbol{\gamma}}$ is the characteristic function of the subinterval $[\mathbf{0}, \boldsymbol{\gamma})$ and, $\lambda_d \big( [\mathbf{0}, \boldsymbol{\gamma}) \big) \coloneqq \prod_{i=1}^{d} \gamma_i$ is the usual $d-$dimensional Lebesque measure. 

\vspace{2mm}
The \textit{$\mathcal{L}_p-$discrepancy} of a sequence $\sigma_N$ is a measure of the irregularity of distribution of $\sigma_N$, and is obtained by taking the $\mathcal{L}_p-$norm $(1 \leq p \leq \infty)$ of the discrepancy function.
\begin{eqnarray}
	\mathcal{L}_{p, N}(\sigma_N) &\coloneqq& \left\Vert g \big( [\mathbf{0}, \boldsymbol{\gamma}), \sigma_N, N \big) \right\Vert _{\mathcal{L}_p} \nonumber \\ &=& \Bigg( \int_{[0, 1)^d} \Big| g \big( [\mathbf{0}, \boldsymbol{\gamma}), \sigma_N, N \big) \Big| ^p d\boldsymbol{\gamma} \Bigg)^{1/p}. \nonumber
\end{eqnarray}

\noindent
Let $\sigma = (\mathbf{b}_n)_{n \in \mathbb{N}} \subset [0, 1)^d$ be an infinite sequence. From the initial segment formed by the first $N$ terms of $\sigma$, we can write $\sigma_N = (\mathbf{b}_i)_{i=1}^N$ and therefore define $\mathcal{L}_{p, N}(\sigma) := \mathcal{L}_{p, N}(\sigma_N)$.

\paragraph{Diaphony.}
In 1976, P. Zinterhof proposed the concept of \textit{diaphony}. It is appropriate that some further notation is now introduced. For any finite sequence $\sigma_N = (\mathbf{a}_i)_{i=1}^N$ contained in $[0, 1)^d$, define the trigonometric sum
\begin{equation*}
S_N(\sigma_N ; \mathbf{m}) \coloneqq \frac{1}{N} \sum_{i=1}^{N} e(\mathbf{m} \cdot \mathbf{a}_i),
\end{equation*}
where we have set $e(x) \coloneqq \text{exp}(2 \pi i x)$ throughout for simplicity. For every lattice point \linebreak $\mathbf{m} = (m_1, ..., m_d) \in \mathbb{Z}^d$, we define $R(\mathbf{m}) \coloneqq \prod_{i=1}^{d} \max \big(1, |m_i| \big).$

\vspace{2mm}
\noindent
Let $\sigma_N$ be a finite sequence contained in $[0, 1)^d$. The diaphony of $\sigma_N$ is defined by,
\begin{equation*}
	F_N(\sigma_N) \coloneqq \Bigg( \sum_{\mathbf{m} \in \mathbb{Z}^d} \frac{\big|S_N(\sigma_N ; \mathbf{m}) \big|^2}{R^2(\mathbf{m})} \Bigg)^{\frac{1}{2}}.
\end{equation*}

\noindent
In the case that $\sigma$ denotes an infinite sequence in $[0, 1)^d$, adopting the same notion as above we truncate $\sigma$ to the finite sequence $\sigma_N$, then set $F_N(\sigma) := F_N(\sigma_N)$.

\paragraph{Dyadic Diaphony.}
The \textit{dyadic diaphony} as introduced in \cite{Hellekalek} is the final measure of irregularity of distribution in which we will be interested. The key difference between the classical diaphony and dyadic diaphony is the replacement of the trigonometric functions with the dyadic Walsh functions.\footnote{It is worth noting that the dyadic diaphony was extended once more to arbitrary bases ($b>2$) using the $b-$adic Walsh function system in \cite{Groz}, named the \textit{$b-$adic diaphony}. See the open problem on page 10.}

\vspace{2mm}
For $k \in \mathbb{N}_0$ with base 2 representation $k = \kappa_{a-1} 2^{a-1} + ... + \kappa_1 2 + \kappa_0$, where $\kappa_i \in \{0, 1\}$ and $\kappa_{a-1} \neq 0$, we define the $k^{th}$ \textit{(dyadic) Walsh function} $\text{wal}_k : \mathbb{R} \rightarrow \{-1, 1\}$, periodic with period one, by
\begin{equation*}
\text{wal}_k(x) \coloneqq (-1)^{x_1 \kappa_0 +... + x_a \kappa_{a-1}},
\end{equation*}
\noindent
for $x \in [0, 1)$ with base 2 representation $x = \frac{x_1}{2} + \frac{x_2}{2^2} + ...$ (unique in the sense that infinitely many of the digits $x_i$ must be zero). For dimension $d \geq 2$, we define the $d-$\textit{dimensional $\mathbf{k}^{th}$ (dyadic) Walsh function} $\text{wal}_{\mathbf{k}} : \mathbb{R}^d \rightarrow \{-1, 1\} $ by
\begin{equation*}
\text{wal}_{\mathbf{k}}(\mathbf{x}) \coloneqq \prod_{j=1}^{d} \text{wal}_{k_j}(x_j),
\end{equation*}
where $\mathbf{k} = (k_1, ..., k_d) \in \mathbb{N}_0^d$ and $\mathbf{x} = (x_1, ..., x_d) \in [0, 1)^d$.  The system $\{\text{wal}_{\mathbf{k}} : \mathbf{k} \in \mathbb{N}_0^d\}$ is called the $d-$dimensional (dyadic) Walsh function system.

The \textit{dyadic diaphony} of a finite sequence $\sigma_N = (\mathbf{a}_i)_{i=1}^{N}$ contained in  $[0, 1)^d$ is defined as,
\begin{equation*}
F_{2, N}(\sigma_N) \coloneqq \Bigg( \frac{1}{3^d-1} \sum_{\mathbf{k} \in \mathbb{N}_0^d \setminus \{\mathbf{0}\}} r_2(\mathbf{k}) \hspace{1mm} \bigg| \frac{1}{N} \sum_{i=1}^{N} \text{wal}_{\mathbf{k}} (\mathbf{a}_i) \bigg|^2 \Bigg)^{1/2}
\end{equation*}
where for $\mathbf{k} = (k_1, ..., k_d) \in \mathbb{N}_0^d$, $r_2(\mathbf{k}) := \prod_{j=1}^{d} r_2(k_j)$, and
\begin{equation*}
r_2(k) \coloneqq  \begin{cases}
1  & \text{if} \hspace{1mm} k=0 \\
2^{-2a} & \text{if} \hspace{1mm} 2^a \leq k < 2^{a+1}, \hspace{1mm} \text{with} \hspace{1mm} a \in \mathbb{N}_0. \\
\end{cases}
\end{equation*}

\noindent
In the scenario that we have an infinite sequence $\sigma \subset [0, 1)^d$, again simply take the initial segment formed by the first $N$ terms of $\sigma$. 

\paragraph{Walsh Series.}
A Walsh system analogue of the trigonometric Fourier series exists, named the \textit{Walsh series} (in some literature, the \textit{Walsh-Fourier Series}). For a function $f: [0, 1)^d \rightarrow \mathbb{R}$, we define the $\mathbf{k}^{th}$ (dyadic) Walsh coefficient of $f$ by,
\begin{equation*}
\widehat{f}(\mathbf{k}) := \int_{[0, 1)^d} f(\mathbf{x}) \text{wal}_{\mathbf{k}} (\mathbf{x}) d\mathbf{x}
\end{equation*}
for $\mathbf{x} \in [0, 1)^d$ and $\mathbf{k} \in \mathbb{N}_0^d$. We can form the Walsh series of $f$ as,
\begin{equation*}
f(\mathbf{x}) \sim \sum_{\mathbf{k} \in \mathbb{N}_0^d} \widehat{f}(\mathbf{k}) \text{wal}_\mathbf{k} (\mathbf{x}).
\end{equation*}

It is appropriate to note that Parseval's identity holds for the Walsh coefficients due to the completeness of the Walsh function system. That is, 
\begin{equation*}
\int_{[0, 1)^d} \big| f(\mathbf{x}) \big|^2 \hspace{1mm} d\mathbf{x} = \sum_{\mathbf{k} \in \mathbb{N}_0^d} \big| \widehat{f}(\mathbf{k}) \big|^2.
\end{equation*}

We refer to [3, Appendix A] for a full treatment of the theory of the Walsh function system and for justification of all above.

\paragraph{Symmetric Sequences.}
Finally, we introduce an important symmetrisation \linebreak technique used in \cite{Proinov 3, Proinov 4}. Let $\sigma_N = (\mathbf{a}_i)_{i=1}^{N}$ be a finite sequence contained in $[0, 1)^d$, and let $\mathbf{x} = (x_1, x_2, ..., x_d) \in [0, 1)^d$. We say that point $\mathbf{x}$ has \textit{multiplicity} $p \hspace{1.5mm} (0 \leq p \leq d)$ with respect to $\sigma_N$, if exactly $p$ terms of $\sigma_N$ coincide with $\mathbf{x}$.

\vspace{1mm}
The sequence $\sigma_N$ is called \textit{symmetric} if for any point $\mathbf{x} = (x_1, ..., x_d) \in [0, 1)^d$, all points of the form,
\begin{equation} \label{sym}
\Big( \tau_1 + (-1)^{\tau_1} x_1, \tau_2 + (-1)^{\tau_2} x_2, ..., \tau_d + (-1)^{\tau_d} x_d \Big) \tag{$**$}
\end{equation}
have the same multiplicity with respect to $\sigma_N$, when $\tau_i \in \{0, 1\}$ independently for $1 \leq i \leq d$. Now let $\tilde{\sigma}_N = (\mathbf{b}_i)_{i=1}^{N}$ be a symmetric sequence contained in $[0, 1)^d$. We say $\tilde{\sigma}_N$ is \textit{generated by} sequence $\sigma_n = (\mathbf{a}_i)_{i=1}^{n}$ if:
\begin{enumerate}
	\item $N = 2^dn$, and
	\item  a point $\mathbf{x} = (x_1, ..., x_d) \in [0, 1)^d$ is a term of the sequence $\sigma_n$, then each point of type (\ref{sym}) is a term of the sequence $\tilde{\sigma}_N$, where $\tau_i \in \{0, 1\}, (1 \leq i \leq d)$ independently.\footnote{Note that every point $x \in [0, 1)^d$ can be regarded as one-term sequence, so every point $x \in [0, 1)^d$ generates at least one symmetric sequence in $[0, 1)^d$ consisting of $p=2^d$ points. Conversely, every symmetric sequence in $[0, 1)^d$ consisting of $p=2^d$ terms is generated by any of its terms.}
\end{enumerate}

See Figure \ref{figure1} and Figure \ref{figure2} below.

Let $\tilde{\sigma} = (\mathbf{b}_n)_{n \in \mathbb{N}}$ be an infinite sequence, $\tilde{\sigma}$ is said to be symmetric if for any $n \in \mathbb{N}$ the finite sequence consisting of $p = 2^d$ terms,
\begin{equation} \label{infsym}
\mathbf{b}_{(n-1)p+1}, \mathbf{b}_{(n-1)p+2}, ..., \mathbf{b}_{np}  \tag{$\dagger$}
\end{equation}
is symmetric. We say that the infinite symmetric sequence $\tilde{\sigma} =  (\mathbf{b}_n)_{n \in \mathbb{N}}$ is generated by an infinite sequence $\sigma = (\mathbf{a}_n)_{n \in \mathbb{N}}$ if for any $n \in \mathbb{N}$, the finite sequence (\ref{infsym}) is generated by the point $\mathbf{a}_n$.

\begin{figure}[H]
	\begin{minipage}{.5\textwidth}
		\centering
		\captionsetup{justification=centering}
		\includegraphics[width=0.9 \linewidth]{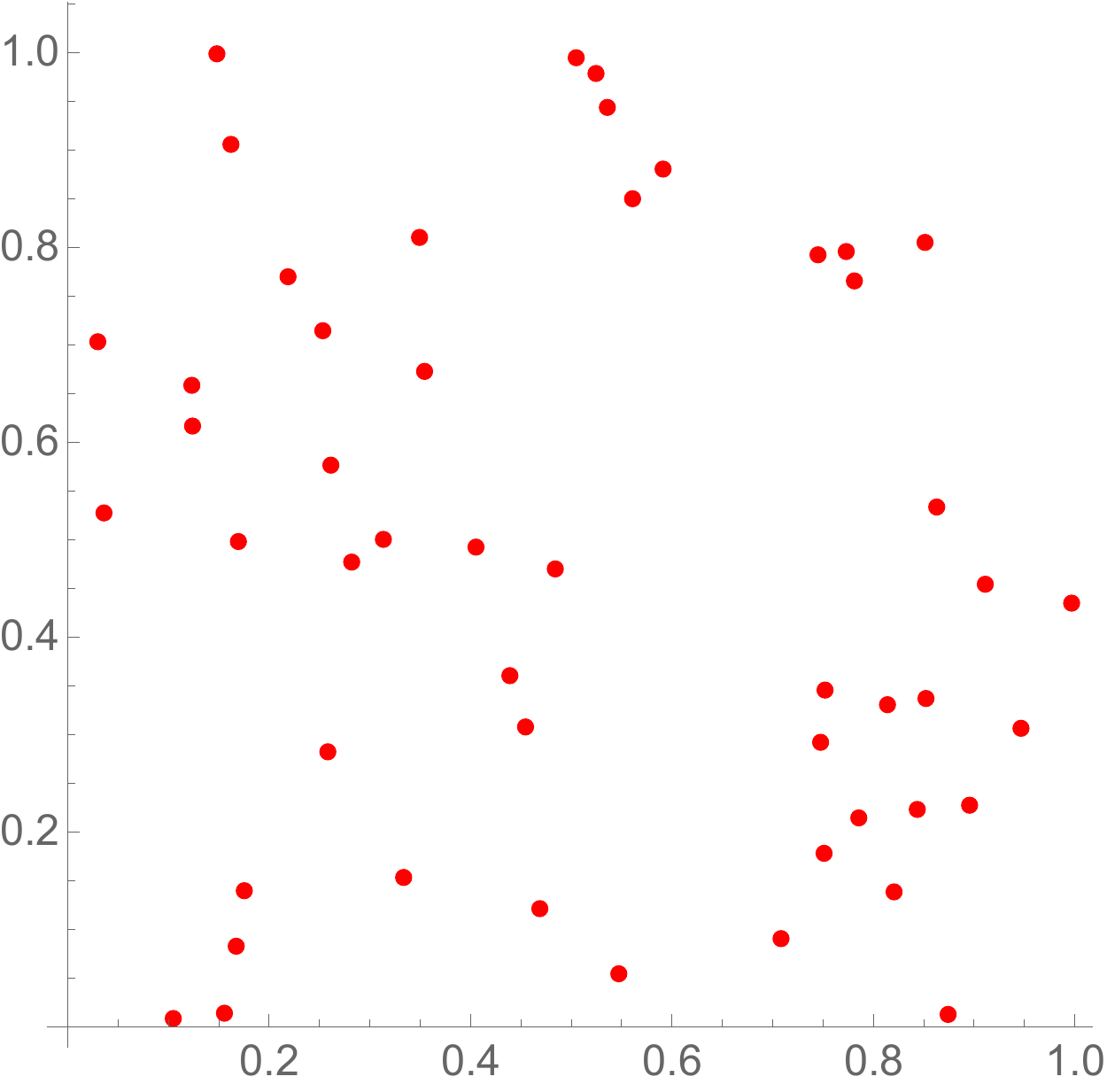}
		\caption{\textit{A random \\ sequence $\sigma_n = (\mathbf{a}_i)_{i=1}^{n} \subset [0, 1)^2$, \\ with $n=50$}}
		\label{figure1}
	\end{minipage}%
	\begin{minipage}{.5\textwidth}
		\centering
		\captionsetup{justification=centering}
		\includegraphics[width=0.9 \linewidth]{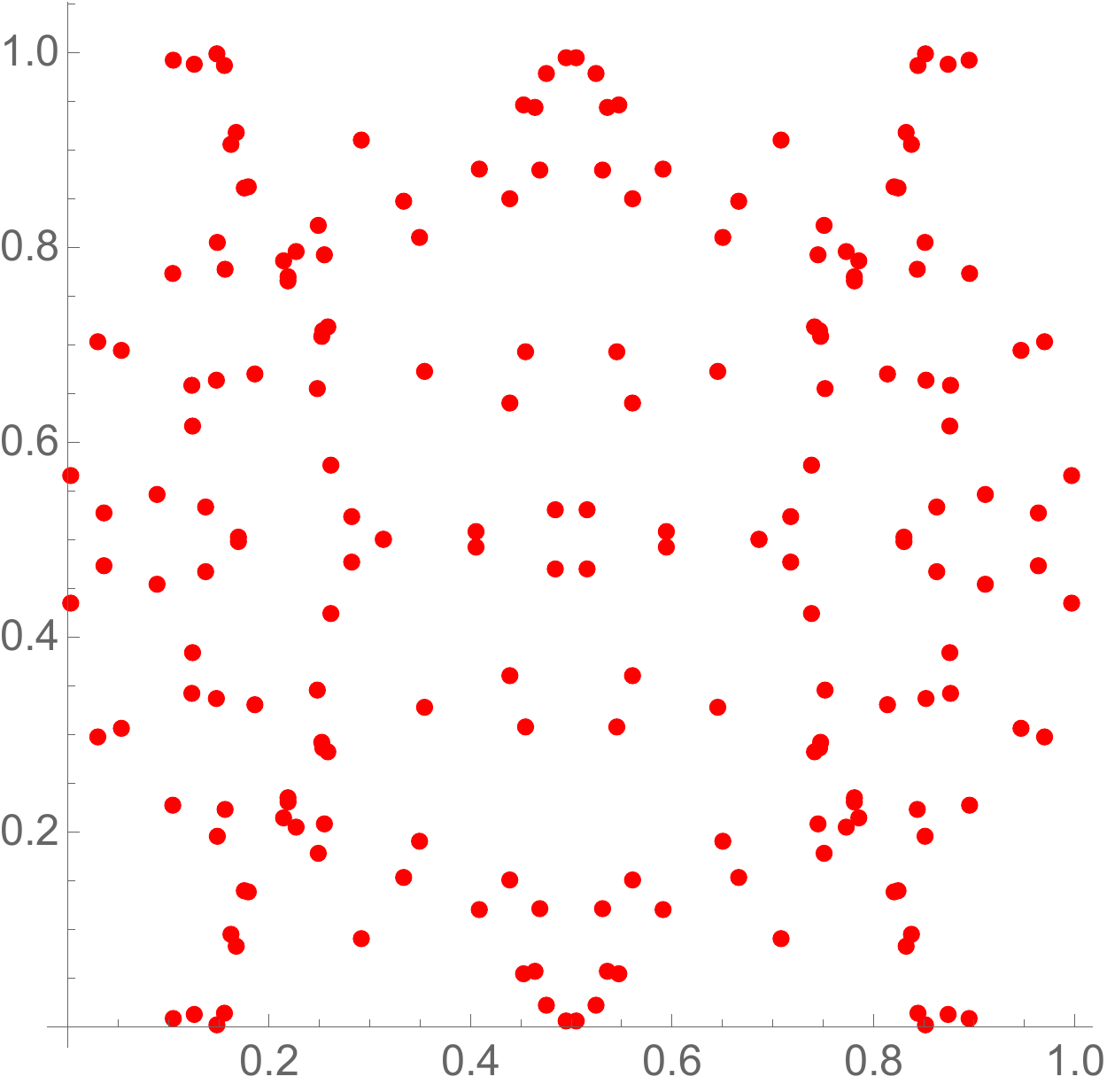}
		\caption{\textit{The symmetric sequence $\tilde{\sigma}_N = (\mathbf{b}_i)_{i=1}^{N} \subset [0, 1)^2$, generated by $\sigma_n$ \\ (from Figure 1)}}
		\label{figure2}
	\end{minipage}
\end{figure}

\begin{remark}
	The above statements regarding generating symmetric sequences have the following equivalent formation.
	
	We say that the symmetric sequence $\tilde{\sigma}_N$ is generated by $\sigma_n = (\mathbf{a}_i)_{i=1}^n$, if every term of $\tilde{\sigma}_N$ can be represented as
	\vspace{1mm}
	\begin{equation*}
	\frac{1}{2}(\mathbf{1} - \boldsymbol{\theta}) + \boldsymbol{\theta}\mathbf{a}_i
	\end{equation*}
	with $1 \leq k \leq n$, $\boldsymbol{\theta} \in Z_d$. $Z_d$ denotes the subset of all $d-$dimensional points of the form $\boldsymbol{\theta} = (\theta_1, ..., \theta_d)$ with each coordinate $\theta_j = \pm 1$ for $1 \leq j \leq d$, and the binary operation between $\boldsymbol{\theta}$ and $\mathbf{a}_i$ is component-wise multiplication.
\end{remark}

\subsection{The Results of Proinov}

Proinov's argument comes in three main steps. An overview of the high-level structure of the proof is contained in Section 3.1. 

\vspace{1mm}
First, Proinov lower bounds the diaphony of a sequence by the $\mathcal{L}_2-$discrepancy of the symmetrised version of the sequence. Theorems 1 and 2 below cater to finite and infinite sequences respectively.

\vspace{1mm}
\begin{theorem}
	Let $\tilde{\sigma}_N$ be any finite symmetric sequence consisting of $N = 2^dn$ terms contained in $[0, 1)^d$, and let $\sigma_n$ be any finite sequence also contained in $[0, 1)^d$ consisting of $n$ terms which generates $\tilde{\sigma}_N$. Then the inequality,
	\begin{equation*}
		\mathcal{L}_{2, N}(\tilde{\sigma}_N) \leq C(d) F_n(\sigma_n)
	\end{equation*}
	holds with
	\vspace{3mm}
	\begin{equation} \label{A1}
	C(d) \coloneqq \frac{1}{2^{d+1}} \sqrt{\bigg( 1-\frac{1}{3^d} \bigg) \bigg( \Big( 1+\frac{6}{\pi^2} \Big)^d-1 \bigg)}.
	\end{equation}
\end{theorem}

\vspace{1mm}
\begin{theorem}
	Let $\tilde{\sigma}$ be any infinite symmetric sequence contained in $[0, 1)^d$, and let $\sigma$ be any infinite sequence contained in $[0, 1)^d$ which generates $\tilde{\sigma}$. Then for a natural number $N \geq 2^d$, the following inequality holds
	\begin{equation*}
		\mathcal{L}_{2, N}(\tilde{\sigma}) \leq C(d) F_n(\sigma) + (2^d-1)/N
	\end{equation*}
	where $n = \lfloor N/2^d \rfloor$ and $C(d)$ defined as in (\ref{A1}).
\end{theorem}

This now allows for the application of the classical lower bound result of Roth. Proinov extends the inequality of Roth to consider infinite sequences contained in $[0, 1)^d$.

\begin{theorem} 
	Let $d \geq 1$. For any infinite sequence $\sigma$ contained in $[0, 1)^d$, we have the following inequality
	\vspace{2mm}
	\begin{equation*}
		\limsup_{N \rightarrow \infty} \frac{N\mathcal{L}_{2, N}(\sigma)}{(\log N)^{\frac{d}{2}}} \geq \alpha(d)
	\end{equation*}
	
	\noindent
	with constant $\alpha(d)$ defined as, 
	\begin{equation} \label{A2}
	\alpha(d) \coloneqq \frac{1}{4^{d+3}(d \log 2)^{\frac{d}{2}}}.
	\end{equation}
\end{theorem}

We take a brief aside at this point to note that in \cite{Pill}, Theorem 3.20 cites a slightly altered constant than $\alpha(d)$ as stated above. In this paper, the author moves forward with constant (\ref{A2}) as defined and used by Proinov to record a self-contained derivation of Proinov's lower bound for the one-dimensional diaphony constant, $f^*$. In any case, the constant $\alpha(d)$ is soon abandoned and replaced by the updated constant $\gamma(d)$ in (\ref{B3}) which is used for the remainder of the text.

\vspace{1mm}
Returning to the results, Proinov combines all the preceding observations to derive his main results regarding the lower bound for the diaphony of finite and infinite sequences in Theorems 4 and 5 respectively.

\begin{theorem} 
	Let $d \geq 2$. For any finite sequence $\sigma_N$ contained in $[0, 1)^d$, we have the following inequality
	\vspace{2mm}
	\begin{equation*}
	\limsup_{N \rightarrow \infty} \frac{N F_N(\sigma_N)}{(\log N)^{\frac{d-1}{2}}} \geq \alpha(d-1)\beta(d)
	\end{equation*}
	with $\alpha(d)$ defined as in (\ref{A2}) and constant $\beta(d)$ defined as,
	\begin{equation} \label{A3}
	\beta(d) \coloneqq 2\pi^d \sqrt{\frac{3^d}{(3^d-1)\big( (\pi^2+6)^d-\pi^{2d}\big)}}.
	\end{equation} 
\end{theorem}

\begin{theorem}
	Let $d \geq 1$. For any infinite sequence $\sigma$ contained in $[0, 1)^d$,
	\begin{equation*}
		\limsup_{N \rightarrow \infty} \frac{N F_N(\sigma)}{(\log N)^{\frac{d}{2}}} \geq \alpha(d)\beta(d)
	\end{equation*}
	where $\alpha(d)$ and $\beta(d)$ are defined in (\ref{A2}) and (\ref{A3}).
\end{theorem}

\subsection{Improvements After 1986}

As a simple Corollary of Theorem 5 (setting $d=1$), the lower bound of the one-dimensional diaphony constant known to Proinov is
\begin{equation} \label{B1}
f^* > \frac{\pi}{256 \sqrt{\log 2}} = 0.0147...
\end{equation}

Looking again at Theorem 5, the constants $\alpha(d)$ and $\beta(d)$ are responsible for arriving at (\ref{B1}). In particular, $\alpha(d)$ originates from the celebrated Theorem of Roth regarding a lower bound for the $\mathcal{L}_2-$discrepancy. The authors in \cite{Aicke} improve this classical result via an adaptation of Roth's method considering certain Fourier coefficients of the discrepancy function with respect to the Haar basis. We formulate this below as Theorem A, however note that the constant $\gamma(d)$ as stated is an edited version to that which was originally published. It was flagged in a later publication \cite{HinLar} by the same co-author that the proof contains a small inaccuracy and instructions are given on how one rectifies this issue, leading to constant (\ref{B3}). For clarity, we attach an appendix which contains the adjusted proof of the $\mathcal{L}_2-$discrepancy result.

\begin{statement}{Theorem A.} [A. Hinrichs \& L. Markhasin, 2011]
Let $d \geq 2$. For a finite sequence $\sigma_N$ contained in $[0, 1)^d$, we have
\begin{equation*} \label{newL2}
\mathcal{L}_{2, N}(\sigma_N) \geq \gamma(d-1) \frac{(\log N)^{\frac{d-1}{2}}}{N}
\end{equation*}
with constant $\gamma(d)$ defined as,
\begin{equation} \label{B3}
\gamma(d) \coloneqq \frac{1}{\sqrt{21} \cdot 2^{2d+1} \sqrt{d!} (\log 2)^{\frac{d}{2}}}.
\end{equation}
\end{statement}

Subsequently, mimicing the proofs of Theorem 3 and Theorem 5 with constant $\alpha(d)$ replaced with $\gamma(d)$, we arrive at the following updated results in the general $d-$dimensional case.

\begin{theorem}
	Let $d \geq 1$. For any infinite sequence $\sigma$ contained in $[0, 1)^d$, we have the following inequality
	\begin{equation*}
	\limsup_{N \rightarrow \infty} \frac{N \mathcal{L}_{2, N}(\sigma)}{(\log N)^{\frac{d}{2}}} \geq \gamma(d)
	\end{equation*}
	with $\gamma(d)$ defined as in (\ref{B3}).
\end{theorem}

\begin{theorem} 
	Let $d \geq 1$. For any infinite sequence $\sigma$ contained in $[0, 1)^d$, we have
	\begin{equation*} 
	\limsup_{N \rightarrow \infty} \frac{N F_N(\sigma)}{(\log N)^{\frac{d}{2}}} \geq \beta(d) \gamma(d)
	\end{equation*}
	with $\beta(d)$ and $\gamma(d)$ are defined as in (\ref{A3}) and (\ref{B3}).
\end{theorem}

One further improvement was made for sequences contained within the unit square, $[0, 1)^2$. The authors in \cite{HinLar} derive an improved lower bound for the $\mathcal{L}_2-$discrepancy of $2-$dimensional finite sequences using a variant of the method from the earlier paper, \cite{Aicke}. For convenience, this is made explicit in Theorem B below and we refer to Section 2 of \cite{HinLar} for a derivation and explicit form of the constant.

\begin{statement}{Theorem B.} [A. Hinrichs \& G. Larcher, 2016]
	For a finite sequence $\sigma_N$ contained in $[0, 1)^2$, the following inequality holds
	\begin{equation*}
	\mathcal{L}_{2, N}(\sigma_N) \geq 0.0515599... \cdot \frac{\sqrt{\log N}}{N}.
	\end{equation*}
\end{statement}

Using a similar argument to that of Theorem 3, it can be shown that one can use the results of $(d+1)-$dimensional finite point sets to study $d-$dimensional infinite sequences. Therefore, we can implement this $2-$dimensional asymptotic constant to gain a most improved lower bound of the one-dimensional diaphony constant.

\begin{corollary}
\begin{equation*}
f^* > 0.0515599... \cdot \pi = 0.1619...
\end{equation*}
\end{corollary}

\subsection{An Extension to the Dyadic Diaphony}

Finally, we apply the technique of Proinov to derive an explicit lower bound for the dyadic diaphony. As above, we consider similar lower bounds for the one-dimensional dyadic diaphony constant which we define as, 
\begin{equation*}
f^*_2 := \inf_{\sigma} \limsup_{N \rightarrow \infty} \frac{NF_{2, N}(\sigma)}{\sqrt{\log N}}.
\end{equation*}

\begin{theorem}
Let $\tilde{\sigma}_N$ be any finite symmetric sequence contained in $[0, 1)^d$ consisting of $N = 2^dn$ terms with $\sigma_n$ any finite sequence contained in $[0, 1)^d$ consisting of $n$ terms which generates $\tilde{\sigma}_N.$ Then,
\begin{equation*}
\mathcal{L}_{2, N}(\tilde{\sigma}_N) \leq \delta(d) F_{2, n}(\sigma_n)
\end{equation*}
holds with constant,
\begin{equation} \label{C1}
\delta(d) \coloneqq \sqrt{3^d-1}.
\end{equation}
\end{theorem}

\begin{theorem}
	Let $\tilde{\sigma}$ be any infinite symmetric sequence contained in $[0, 1)^d$, and let $\sigma$ be any infinite sequence contained in $[0, 1)^d$ which generates $\tilde{\sigma}$. Then for a natural number $N \geq 2^d$,
	\begin{center}
		$\mathcal{L}_{2, N}(\tilde{\sigma}) \leq \delta(d) F_{2, n}(\sigma) + (2^d-1)/N$
	\end{center}
	where $n = \lfloor N/2^d \rfloor$ and $\delta(d)$ defined as in (\ref{C1}).
\end{theorem}

\begin{theorem}
	Let $d \geq 2$. For any finite sequence $\sigma_N$ contained in $[0, 1)^d$, the following inequality holds
	\begin{equation*}
		\limsup_{N \rightarrow \infty} \frac{N F_{2, N}(\sigma_N)}{(\log N)^{\frac{d-1}{2}}} \geq \gamma(d-1) \mu(d)
	\end{equation*}
	with $\gamma(d)$ defined as in (\ref{B3}) and
	\begin{equation} \label{C2}
	\mu(d) \coloneqq \frac{1}{2^d\sqrt{3^d-1}}.
	\end{equation}
\end{theorem}

\begin{theorem}
	Let $d \geq 1$. For any infinite sequence $\sigma$ contained in $[0, 1)^d$,
	\begin{equation*}
		\limsup_{N \rightarrow \infty} \frac{N F_{2, N}(\sigma)}{(\log N)^{\frac{d}{2}}} \geq \gamma(d) \mu(d)
	\end{equation*}
	with $\gamma(d)$ and $\mu(d)$ defined as in (\ref{B3}) and (\ref{C2}).
\end{theorem}

Theorems 10 and 11 allow the computation of exact lower bounds for general $d-$dimensional sequences. However, once again we can use Theorem B  to calculate the best lower bound for the newly defined one-dimensional dyadic diaphony constant.

\begin{corollary}
\begin{equation*}
f^*_2 > 0.0515599... \cdot \frac{1}{2 \sqrt{2}} = 0.0182...
\end{equation*}
\end{corollary}

\vspace{3mm}
\noindent
\textbf{\textit{Open Problem}}
\vspace{1mm}

We leave it as an open problem to derive similar lower bounds using Proinov's methods for the $b-$adic diaphony of sequences contained in the $d-$dimensional unit cube. (See Footnote 3 in Section 2.1).

\vspace{2mm}
As a first step, one would need to formulate a similar inequality to those of Theorems 1 and 8 giving a relationship between the $b-$adic diaphony and the $\mathcal{L}_2-$discrepancy, explicitly forming a constant similar to $C(d)$ or $\delta(d)$ respectively. It is reasonable to conjecture that this constant would depend on (and only on) the dimension $d$, and the choice of base $b$.

\section{The Proofs of Proinov}

In this Section, we present the proofs of Proinov. To the best of our knowledge, with the exception of Theorem 3, the only record of these proofs are contained in Proinov's monograph \cite{Proinov 2} which is written in Bulgarian and not widely available. Proinov's proof of Theorem 3 is given in full as Theorem 2.2 in \cite{DickPill}, we refer the curious reader to this survey.

\vspace{1mm}
In Section 3.1, we outline the argument of Proinov since his method is of general interest and will also be used in Section 4 to derive a lower bound for the dyadic diaphony. Section 3.2 contains the full proofs of Theorems 1, 2, 4 and 5.
 
\subsection{The Main Ideas of Proinov}

We outline the major steps used to formulate Theorem 4, the Theorem concerning explicit lower bounds for the diaphony of finite sequences contained in $[0, 1)^d$. The extension to infinite sequences (to derive Theorem 5) follows from several technical Lemmas, the details of which are outlined in the next subsection. 

\vspace{1mm}
In the first instance, Proinov formulates Theorem 1 which lower bounds the diaphony of a sequence by the $\mathcal{L}_2-$discrepancy of the symmetrised version of the sequence. That is, for finite sequences $\sigma_n$ and symmetric $\tilde{\sigma}_N$ contained in $[0, 1)^d$ consisting of $n$ and $N = 2^dn$ terms respectively such that $\sigma_n$ generates $\tilde{\sigma}_N$, we have
\begin{equation*}
\mathcal{L}_{2, N}(\tilde{\sigma}_N) \leq C(d) F_n(\sigma_n)
\end{equation*}
holds with constant $C(d)$ defined as previously in (\ref{A1}). There is significant machinery involved in deducing this result. Specifically, the discrepancy function $g \big( [\mathbf{0}, \boldsymbol{\gamma}), \tilde{\sigma}_N, N \big)$ defined in the introduction is expanded as a Fourier series,
\begin{equation*}
g \big( [\mathbf{0}, \boldsymbol{\gamma}), \tilde{\sigma}_N, N \big) \sim \sum_{\mathbf{m} \in \mathbb{Z}^d} \widehat{g}(\mathbf{m}) e(\mathbf{m} \cdot \boldsymbol{\gamma})
\end{equation*}
where $\widehat{g}(\mathbf{m})$ denote the Fourier coefficients. Proinov then implements Parseval's identity with the discrepancy function to obtain an expression for the $\mathcal{L}_2-$discrepancy in terms of the Fourier coefficients,
\begin{equation*}
\int_{[0, 1)^d} \Big| g \big( [\mathbf{0}, \boldsymbol{\gamma}), \tilde{\sigma}_N, N \big) \Big| ^2 d\boldsymbol{\gamma} = \mathcal{L}_{2,N}^2(\tilde{\sigma}_N) = \sum_{\mathbf{m} \in \mathbb{Z}^d} \big| \widehat{g}(\mathbf{m}) \big | ^2.
\end{equation*}

\noindent
Subsequently, with some rigorous calculation one finds that the summation above can be approximated as,
\begin{equation*}
\sum_{\mathbf{m} \in \mathbb{Z}^d} \big| \widehat{g}(\mathbf{m}) \big | ^2 \leq C^2(d) \sum_{\mathbf{m} \in \mathbb{Z}^d} \frac{\big|S_n(\sigma_n ; \mathbf{m}) \big|^2}{R^2(\mathbf{m})}
\end{equation*}
with all terms on the right-hand side of this inequality defined as in Section 2.1, and $C(d)$ as in (\ref{A1}). Putting the last two lines together obtains the desired result for Theorem 1,
\begin{equation*}
\mathcal{L}_{2, N}^2(\tilde{\sigma}_N) \leq C^2(d) F^2_n(\sigma_n).
\end{equation*}

\noindent
Rearranging Theorem 1, taking Roth's 1954 classical lower bound for the $\mathcal{L}_2-$discrepancy while noting that
\begin{equation*}
\frac{1}{C(d)} = 2^d\beta(d),
\end{equation*}
where $\beta(d)$ is defined as (\ref{A3}), we arrive with some small manipulation at
\begin{equation*}
F_n(\sigma_n) \geq \alpha(d-1)\beta(d)\frac{(\log n)^{\frac{d-1}{2}}}{n}
\end{equation*}
as required with an explicit lower bound for the diaphony of an arbitrary finite sequence contained in $[0, 1)^d$.

\subsection{The Proofs of Theorems 1, 2, 4 and 5}

\paragraph{Proof of Theorem 1.} Let $\sigma_n = (\mathbf{a}_i)_{i=1}^n \subset [0, 1)^d$ be any finite sequence which generates a symmetric sequence $\tilde{\sigma}_N = (\mathbf{b}_i)_{i=1}^N \subset [0, 1)^d$, containing $n$ and $N=2^dn$ terms respectively. From the definition of the sequence $\tilde{\sigma}_N$ which is generated by $\sigma_n$,
\begin{equation*}
\frac{1}{N} \sum_{i=1}^{N} \chi_{\boldsymbol{\gamma}} (\mathbf{b}_i) = \frac{1}{2^dn} \sum_{i=1}^{n} \sum_{\boldsymbol{\theta} \in Z_d} \chi_{\boldsymbol{\gamma}} \bigg( \frac{1}{2}(\mathbf{1} - \boldsymbol{\theta}) + \boldsymbol{\theta} \mathbf{a}_i \bigg) \nonumber
\end{equation*}
where the set $Z_d$ is defined as in Remark 1. Therefore, we rewrite the discrepancy function as
\begin{equation} \label{eq:5}
g \big( [\mathbf{0}, \boldsymbol{\gamma}), \tilde{\sigma}_N, N \big) = \frac{1}{2^dn} \sum_{i=1}^{n} \sum_{\boldsymbol{\theta} \in Z_d} \chi_{\boldsymbol{\gamma}} \bigg( \frac{1}{2}(\mathbf{1} - \boldsymbol{\theta}) + \boldsymbol{\theta} \mathbf{a}_i \bigg) - \lambda_d \big( [\mathbf{0}, \boldsymbol{\gamma}) \big).
\end{equation}

\noindent
The function can be written asymptotically equal to a Fourier series,
\begin{equation*} 
g \big( [\mathbf{0}, \boldsymbol{\gamma}), \tilde{\sigma}_N, N \big) \sim \sum_{\mathbf{m} \in \mathbb{Z}^d} \widehat{g}(\mathbf{m}) e(\mathbf{m} \cdot \boldsymbol{\gamma})
\end{equation*}
where $\widehat{g}(\mathbf{m})$ denote the Fourier coefficients. Each $\widehat{g}(\mathbf{m})$ can be calculated by,
\begin{equation} \label{eq:6}
\widehat{g}(\mathbf{m}) = \int_{[0, 1)^d} g \big( [\mathbf{0}, \boldsymbol{\gamma}), \tilde{\sigma}_N, N \big) e(\mathbf{-m} \cdot \boldsymbol{\gamma}) d\boldsymbol{\gamma}.
\end{equation}

\noindent
We need the following one-dimensional integrals defined and denoted,
\begin{equation*}
A(m) := \int_{0}^{1} \gamma e(-m\gamma)d\gamma \hspace{5mm} \& \hspace{5mm} B(m, a) := \int_{a}^{1} e(-m\gamma)d\gamma
\end{equation*}
for $m \in \mathbb{Z}$ and $a \in E$. It is easily calculated that,
\begin{equation} \label{eq:7}
A(m) = \begin{cases}
-\frac{1}{2\pi i m}, & \text{if} \hspace{1mm} m \neq 0 \\
\frac{1}{2}, & \text{if} \hspace{1mm} m = 0
\end{cases}
\end{equation}
\begin{center}
	\&
\end{center} 
\begin{equation} \label{eq:8}
B(m, a) = \begin{cases}
\big( \frac{1}{2\pi i m} \big) \big( e(-ma) - 1 \big), & \text{if} \hspace{1mm} m \neq 0 \\
1 - a, & \text{if} \hspace{1mm} m = 0.
\end{cases}
\end{equation}

\noindent
Using equations (\ref{eq:5}) to (\ref{eq:8}) and by writing
\begin{equation*}
\boldsymbol{\zeta}_{i}(\boldsymbol{\theta}) := \frac{1}{2} (\mathbf{1} - \boldsymbol{\theta}) + \boldsymbol{\theta} \mathbf{a}_{i},
\end{equation*}
we obtain the following expression for the Fourier coefficients.
\begin{eqnarray} \label{eq:9}
\widehat{g}(\mathbf{m}) &=&  \frac{1}{2^dn} \sum_{i=1}^{n} \sum_{\boldsymbol{\theta} \in Z_d} \int_{[0, 1)^d} \chi_{\boldsymbol{\gamma}}(\boldsymbol{\zeta}_{i}) e(-\mathbf{m} \cdot \boldsymbol{\gamma}) d \boldsymbol{\gamma} - \int_{[0, 1)^d} \lambda_d\big( [ \mathbf{0}, \boldsymbol{\gamma}) \big) e(-\mathbf{m} \cdot \boldsymbol{\gamma}) d\boldsymbol{\gamma} \nonumber \\ &=&  \frac{1}{2^dn} \sum_{i=1}^{n} \sum_{\boldsymbol{\theta} \in Z_d} \int_{\boldsymbol{\zeta_{i}}}^{\mathbf{1}} e(-\mathbf{m} \cdot \boldsymbol{\gamma}) d\boldsymbol{\gamma} - \int_{[0, 1)^d} \lambda_d\big( [\mathbf{0}, \boldsymbol{\gamma}) \big) e(-\mathbf{m} \cdot \boldsymbol{\gamma}) d\boldsymbol{\gamma} \nonumber \\ &=& \frac{1}{2^dn} \sum_{i=1}^{n} \sum_{\boldsymbol{\theta} \in Z_d} \prod_{j=1}^{d} \int_{\zeta_{ij}}^{1} e(-m_j \gamma_j) d\gamma_j - \prod_{j=1}^{d} \int_{0}^{1} \gamma_j \hspace{0.5mm} e(-m_j \gamma_j) d\gamma_j \nonumber \\ &=& \frac{1}{2^dn} \sum_{i=1}^{n}  \prod_{j=1}^{d} \sum_{\theta_j = \pm 1} B(m_j , \zeta_{ij}) - \prod_{j=1}^{d} A(m_j).
\end{eqnarray}

\noindent
Note that
\begin{equation*}
\widehat{g}(\mathbf{0}) = \frac{1}{2^d} - \frac{1}{2^dn} \sum_{i=1}^{n} 1 = 0,
\end{equation*}
which follows from 
\begin{eqnarray} \label{eq:10}
\sum_{\theta_j = \pm 1} B(0, \zeta_{ij}) &=& \sum_{\theta_j = \pm 1} 1 - \zeta_{ij} \nonumber \\ &=& \sum_{\theta_j = \pm 1} 1 - \Bigg( \frac{1}{2}(1 - \theta_j) + \theta_j a_{ij} \Bigg) \nonumber \\ &=& ( 1 - a_{ij} ) + (1 - 1 + a_{ij} ) = 1.
\end{eqnarray}

\noindent
Using Parseval's identity,
\begin{eqnarray} \label{eq:11}
\mathcal{L}_{2,N}^2(\tilde{\sigma}_N) &=& \int_{[0, 1)^d} \Big| g \big( [\mathbf{0}, \boldsymbol{\gamma}), \tilde{\sigma}_N, N \big) \Big| ^2 d\boldsymbol{\gamma} \nonumber \\ &=& \sum_{\mathbf{m} \in \mathbb{Z}^d} \big| \widehat{g}(\mathbf{m}) \big | ^2 = \sideset{}{'}\sum_{\mathbf{m} \in \mathbb{Z}^d} \big|  \widehat{g}(\mathbf{m}) \big| ^2
\end{eqnarray} 
where $\sideset{}{'} \sum$ denotes the sum without the zero index.

\vspace{1mm}
\noindent
Let $A$ be an arbitrary nonempty subset of $\{1, 2, ..., d\}$. Denote by $M(A)$, the set consisting of integer points $\mathbf{m} = (m_1, ..., m_d)$ such that $m_j \neq 0 \hspace{1mm} (1 \leq j \leq d)$ if and only if $j \in A$. We define,
\begin{equation} \label{eq:12}
\phi(A) := \sum_{\mathbf{m} \in M(A)} \big| \widehat{g}(\mathbf{m}) \big| ^2
\end{equation}

\noindent
and it is therefore easy to see that (\ref{eq:11}) can be written in the following form,
\begin{eqnarray} \label{eq:13}
\mathcal{L}_{2,N}^2(\tilde{\sigma}_N) &=& \sum_{p=1}^{d} \sum_{|A| = p} \sum_{\mathbf{m} \in M(A)} \big| \widehat{g}(\mathbf{m}) \big| ^2 \nonumber \\ &=& \sum_{p=1}^{d} \sum_{|A| = p} \phi(A).
\end{eqnarray}
The sum is over all subsets $A$ of $\{1, 2, ...., d\}$ such that $|A| = p$, for $1 \leq p \leq d$.

\vspace{1mm}
\noindent
Now fix $A$ as some nonempty subset of $\{1, 2, ..., d\}$ with $p$ elements. We prove the estimate, 
\begin{equation} \label{eq:14}
\phi(A) \leq C^2(d) \sum_{\mathbf{m} \in M(A)} \frac{\big| S_n(\sigma_n ; \mathbf{m}) \big|^2}{R^{2}(\mathbf{m})} 
\end{equation}
with constant $C(d)$ as in (\ref{A1}). Returning to the expression (\ref{eq:9}) for the Fourier coefficients, take $\mathbf{m} \in M(A)$. It follows from the work done so far that,
\begin{equation} \label{eq:15}
\widehat{g}(\mathbf{m}) = \frac{1}{2^dn(2 \pi i )^p R(\mathbf{m})} \sum_{i=1}^{n} \prod_{j \in A} \sum_{\theta_j = \pm 1} \big( e(- \theta_j m_j a_{ij}) - 1 \big) - \frac{(-1)^p}{(2\pi i)^p 2^{d-p} R(\mathbf{m})}
\end{equation}

\noindent
and we make the following transformation,
\begin{eqnarray} \label{eq:16}
\prod_{j \in A} \sum_{\theta_j = \pm 1} \big( e(- \theta_j m_j a_{ij}) - 1 \big) &=& \prod_{j \in A} \big( e(m_j a_{ij}) + e(- m_j a_{ij})- 2 \big) \nonumber \\ &=& \sum_{\boldsymbol{\epsilon} \in E(A)} r(\boldsymbol{\epsilon}) e(\boldsymbol{\epsilon} \mathbf{m}  \cdot \mathbf{a}_i).
\end{eqnarray}

\noindent
$r$ is a coefficient function and $E(A)$ denotes all points $\boldsymbol{\epsilon} = (\epsilon_1, ..., \epsilon_d)$ such that $\epsilon_j \hspace{1mm} (1 \leq j \leq d)$ is equal to one of $-1, 0, 1$ if $j \in A$ and $\epsilon_j = 0$ otherwise. The operation between $\boldsymbol{\epsilon}$ and $\mathbf{m}$ is component-wise multiplication.
\noindent
Note that clearly $|E(A)| = 3^p$ since $|A| = p$, and moreover it is easily seen that
\begin{equation} \label{eq:17}
r(\mathbf{0}) = (-2)^p
\end{equation} 
and
\begin{equation} \label{eq:18}
| r(\boldsymbol{\epsilon}) | \leq 2^{p-1}
\end{equation}
for all $\boldsymbol{\epsilon} \neq \mathbf{0}$. Now from (\ref{eq:15}), (\ref{eq:16}) and (\ref{eq:17}),
\begin{eqnarray}
\widehat{g}(\mathbf{m}) &=& \frac{1}{2^dn(2 \pi i)^p R(\mathbf{m})} \sum_{i=1}^{n} \sum_{\boldsymbol{\epsilon} \in E(A)} r(\boldsymbol{\epsilon}) e(\boldsymbol{\epsilon} \mathbf{m} \cdot \mathbf{a}_i) - \frac{(-1)^p}{(2\pi i)^p 2^{d-p} R(\mathbf{m})} \nonumber \\ &=& \frac{1}{2^d(2 \pi i)^p R(\mathbf{m})} \sideset{}{'}\sum_{\boldsymbol{\epsilon} \in E(A)} r(\boldsymbol{\epsilon}) S_n(\sigma_n ; \boldsymbol{\epsilon} \mathbf{m}) + \frac{(-2)^p}{(2\pi i)^p 2^d R(\mathbf{m})} - \frac{(-1)^p}{(2\pi i)^p 2^{d-p} R(\mathbf{m})} \nonumber \\ &=& \frac{1}{2^d(2 \pi i)^p R(\mathbf{m})} \sideset{}{'}\sum_{\boldsymbol{\epsilon} \in E(A)} r(\boldsymbol{\epsilon}) S_n(\sigma_n ; \boldsymbol{\epsilon} \mathbf{m}) \nonumber.
\end{eqnarray}

\noindent
By this last equality and (\ref{eq:18}), we obtain the following estimate for the Fourier coefficients of the discrepancy function
\begin{eqnarray}
\big| \widehat{g}(\mathbf{m}) \big | &\leq& \frac{2^{p-1}}{2^d(2\pi)^p} \Bigg| \sideset{}{'}\sum_{\boldsymbol{\epsilon} \in E(A)}  \frac{S_n(\sigma_n ; \boldsymbol{\epsilon} \mathbf{m})}{R(\mathbf{m})} \Bigg| \nonumber \\ &\leq& \frac{1}{2^{d+1} \pi^p} \sideset{}{'}\sum_{\boldsymbol{\epsilon} \in E(A)} \frac{\big| S_n(\sigma_n ; \boldsymbol{\epsilon} \mathbf{m}) \big|}{R(\mathbf{m})} \nonumber,
\end{eqnarray}

\noindent
and using the Cauchy-Schwarz inequality on the right-hand side of the above gives,
\begin{equation*}
\big| \widehat{g}(\mathbf{m}) \big| ^2 \leq \frac{3^p - 1}{4^{d+1} \pi^{2p}} \sideset{}{'}\sum_{\boldsymbol{\epsilon} \in E(A)}  \frac{\big|S_n(\sigma_n ; \boldsymbol{\epsilon} \mathbf{m}) \big| ^2}{R^{2}(\mathbf{m})}
\end{equation*}
recalling that $|E(A)| = 3^p$. Returning to (\ref{eq:12}),
\begin{equation} \label{eq:19}
\phi(A) \leq \frac{3^p - 1}{4^{d+1} \pi^{2p}} \sideset{}{'}\sum_{\boldsymbol{\epsilon} \in E(A)} \Omega(\boldsymbol{\epsilon})
\end{equation} 
where we have written
\begin{equation*}
\Omega(\boldsymbol{\epsilon}) := \sum_{\mathbf{m} \in M(A)} \frac{\big|S_n(\sigma_n ; \boldsymbol{\epsilon} \mathbf{m}) \big| ^2}{R^{2}(\mathbf{m})}.
\end{equation*}

\noindent
For any nonempty subset $B$ of the set $\{1, 2, ..., d\}.$ We introduce the set $U(B)$, consisting of all points $\boldsymbol{\epsilon} = (\epsilon_1, ..., \epsilon_d)$ such that $\epsilon_j = \pm1$ if $j \in B$ and $\epsilon_j = 0$ otherwise. Clearly, $|U(B)| = 2^{|B|}.$ We can now identify
\begin{equation} \label{eq:20}
\sideset{}{'}\sum_{\boldsymbol{\epsilon} \in E(A)} \Omega(\boldsymbol{\epsilon}) = \sum_{q=1}^{p} \sum_{\substack{B \subset A \\ |B| = q}} \sum_{\boldsymbol{\epsilon} \in U(B)} \Omega(\boldsymbol{\epsilon})
\end{equation}
where the summation on the right hand side is over all possible subsets $B$ of $A$ consisting of $q$ elements ($1 \leq q \leq p$).

\noindent
Now let $B$ be a fixed nonempty $q$ element subset of $A$. Then for $\boldsymbol{\epsilon} \in U(B)$,
\begin{eqnarray} \label{eq:21}
\Omega(\boldsymbol{\epsilon}) &=& \sum_{\mathbf{m} \in M(A)} \frac{\big|S_n(\sigma_n ; \boldsymbol{\epsilon} \mathbf{m}) \big| ^2}{R^{2}(\mathbf{m})} \nonumber \\ &=& \Bigg( \sideset{}{'}\sum_{m=- \infty}^{\infty} \frac{1}{m^2} \Bigg)^{p-q} \sum_{\mathbf{m} \in M(B)}  \frac{\big|S_n(\sigma_n ; \boldsymbol{\epsilon} \mathbf{m}) \big| ^2}{R^{2}(\mathbf{m})} \nonumber \\ &=& \bigg( \frac{\pi^2}{3}\bigg)^{p-q} \sum_{\mathbf{m} \in M(B)}  \frac{\big|S_n(\sigma_n ; \boldsymbol{\epsilon} \mathbf{m}) \big| ^2}{R^{2}(\mathbf{m})} \nonumber \\ &=& \bigg(\frac{\pi^2}{3}\bigg)^{p-q} \sum_{\mathbf{m} \in M(B)} \frac{\big|S_n(\sigma_n ; \mathbf{m}) \big| ^2}{R^{2}(\mathbf{m})}
\end{eqnarray}
where the last equality holds since $\boldsymbol{\epsilon} \in U(B)$ has the effect of permuting the elements of the set $M(B)$ in the summation. Therefore from (\ref{eq:20}) and (\ref{eq:21}), we have
\begin{equation} \label{eq:22}
\sideset{}{'}\sum_{\boldsymbol{\epsilon} \in E(A)} \Omega(\boldsymbol{\epsilon}) \leq \frac{\pi^{2p} D(p)}{3^p} \sum_{\mathbf{m} \in M(A)} \frac{\big|S_n(\sigma_n ; \mathbf{m}) \big|^2}{R^2(\mathbf{m})},
\end{equation}
where
\begin{eqnarray} \label{eq:23}
D(p) &=& \sum_{q=1}^{p} \bigg( \frac{3}{\pi^2} \bigg)^q \sum_{\substack{B \subset A \\ |B| = q}} \sum_{\boldsymbol{\epsilon} \in U(B)} 1 \nonumber \\ &=&  \sum_{q=1}^{p} \bigg( \frac{6}{\pi^2} \bigg)^q \sum_{\substack{B \subset A \\ |B| = q}} 1 \nonumber \\ &=& \sum_{q=1}^{p} \binom{p}{q} \bigg( \frac{6}{\pi^2} \bigg)^q = \bigg(1+ \frac{6}{\pi^2} \bigg)^p - 1.
\end{eqnarray}
From (\ref{eq:19}), (\ref{eq:22}) and (\ref{eq:23}), conclude that
\begin{equation*}
\phi(A) \leq \frac{1}{4^{d+1}} \bigg( 1 - \frac{1}{3^p} \bigg) \bigg( \Big( 1+\frac{6}{\pi^2} \Big)^p - 1 \bigg) \sum_{\mathbf{m} \in M(A)} \frac{\big|S_n(\sigma_n ; \mathbf{m}) \big|^2}{R^2(\mathbf{m})}.
\end{equation*}
Hence bearing in mind that for $p \leq d$,
\begin{equation*}
\frac{1}{4^{d+1}} \bigg( 1 - \frac{1}{3^p} \bigg) \bigg( \Big( 1+\frac{6}{\pi^2} \Big)^p - 1 \bigg) \leq C^2(d),
\end{equation*}
and the assertion (\ref{eq:14}) is proved.
\noindent
Now, to finish from (\ref{eq:13}) and (\ref{eq:14})
\begin{eqnarray}
\mathcal{L}_{2, N}^2(\tilde{\sigma}_N) &\leq& C^2(d) \sum_{p=1}^{d} \sum_{|A| = p} \sum_{\mathbf{m} \in M(A)} \frac{\big|S_n(\sigma_n ; \mathbf{m}) \big|^2}{R^2(\mathbf{m})} \nonumber \\ &=& C^2(d) \sum_{\mathbf{m} \in \mathbb{Z}^d} \frac{\big|S_n(\sigma_n ; \mathbf{m}) \big|^2}{R^2(\mathbf{m})} \nonumber \\ &=& C^2(d) F^2_n(\sigma_n) \nonumber
\end{eqnarray}  
and by square rooting, we have the statement. \hspace*{\fill}\qed

\noindent
\paragraph{Proof of Theorem 2.} Let $\tilde{\sigma}$ be an infinite symmetrical sequence contained in $[0, 1)^d$, and let $\sigma$ be an infinite sequence contained in $[0, 1)^d$ which generates $\tilde{\sigma}$. Let a natural number $N \geq 2^d$, then set $n = \big\lfloor N/2^d \big\rfloor$ and $m = 2^dn$. Firstly, notice that
\begin{equation} \label{eq:24}
2^dn \leq N < 2^d(n+1).
\end{equation}
From the definitions of symmetrisation in Section 2.1, set $\tilde{\sigma}_m$ to be the sequence
\begin{center}
	$\tilde{\sigma}_m = (\mathbf{b}_1, \mathbf{b}_2, ..., \mathbf{b}_m)$,
\end{center}
consisting of the first $n$ terms of sequence $\sigma$. Applying Theorem 1 to sequences $\tilde{\sigma}_m$ and $\sigma_n$, we have
\begin{equation} \label{eq:25}
\mathcal{L}_{2, m}(\tilde{\sigma}) = \mathcal{L}_{2, m}(\tilde{\sigma}_m) \leq C(d)F_n(\sigma_n) = C(d)F_n(\sigma),
\end{equation}
where $C(d)$ is as in (\ref{A1}).

\vspace{2mm}
The next step requires a well known technical Lemma.
\begin{lemma} \label{lem:2}
	Let $1 \leq p \leq \infty.$ Let $\sigma$ be any infinite sequence contained in $[0, 1)^d$. Then for every $n \in \mathbb{N}$ with $n \leq N$, the following inequality holds,
	\begin{center}
		$N \mathcal{L}_{p, N}(\sigma) \leq n \mathcal{L}_{p, n}(\sigma) + N - n$.
	\end{center}
\end{lemma}

	\noindent
	\textit{Proof.}	Let $n \in \mathbb{N}$ such that $1 \leq n \leq N$, and $\mathbf{x} \in [0, 1)^d$ be an arbitrary point. Write $\sigma = (\mathbf{a}_n)_{n \in \mathbb{N}}$, then
	\begin{equation} \label{eq:26}
	\sum_{i=1}^{N} \chi_{\mathbf{x}}(\mathbf{a}_i) = \sum_{i=1}^{n} \chi_{\mathbf{x}}(\mathbf{a}_i) + p(\mathbf{x})
	\end{equation}
	where the function $p$ satisfies the condition 
	\begin{equation}  \label{eq:27}
	0 \leq p(\mathbf{x}) \leq N-n.
	\end{equation}
	Using (\ref{eq:26}), the discrepancy function can be written as
	\begin{eqnarray} \label{eq:28}
	Ng\big([\mathbf{0}, \mathbf{x}), \sigma, N \big) &=& \sum_{i=1}^{N} \chi_{\mathbf{x}}(\mathbf{a}_i) - N\lambda_d \big( [\mathbf{0}, \mathbf{x}) \big) \nonumber \\ &=& \sum_{i=1}^{n} \chi_{\mathbf{x}}(\mathbf{a}_i) + p(\mathbf{x}) - N \lambda_d \big( [\mathbf{0}, \mathbf{x}) \big) \nonumber \\ &=& \sum_{i=1}^{n} \chi_{\mathbf{x}}(\mathbf{a}_i) - n\lambda_d \big( [\mathbf{0}, \mathbf{x}) \big) + n\lambda_d \big( [\mathbf{0}, \mathbf{x}) \big) + p(\mathbf{x}) - N\lambda_d \big( [\mathbf{0}, \mathbf{x}) \big) \nonumber \\ &=& \Bigg[ \sum_{i=1}^{n} \chi_{\mathbf{x}}(\mathbf{a}_i) - n\lambda_d \big( [\mathbf{0}, \mathbf{x}) \big) \Bigg] + q(\mathbf{x}) \nonumber \\ &=& n g \big( [\mathbf{0}, \mathbf{x}), \sigma, n \big) + q(\mathbf{x}),
	\end{eqnarray}
	where the function $q$ denotes, 
	\begin{equation}  \label{eq:29}
	q(\mathbf{x}) := p(\mathbf{x}) - (N-n)\lambda_d \big( [\mathbf{0}, \mathbf{x}) \big).
	\end{equation}
	From (\ref{eq:27}), (\ref{eq:29}) and noting that, 
	\begin{center}
		$0 \leq \lambda_d \big( [\mathbf{0}, \mathbf{x}) \big) \leq 1$
	\end{center}
	we can conclude 
	\begin{equation} \label{eq:30}
	\big| q(\mathbf{x}) \big| \leq N-n.
	\end{equation}
	Now we use (\ref{eq:28}), (\ref{eq:30}) and the definition of the $\mathcal{L}_p-$discrepancy to imply
	\begin{eqnarray}
	N \mathcal{L}_{p, N}(\sigma) &=& N\left\Vert g( \cdot , \sigma, N) \right\Vert _{\mathcal{L}_p} \nonumber \\ &=& \left\Vert n g( \cdot , \sigma, n) + q(\cdot) \right\Vert _{\mathcal{L}_p} \nonumber \\ &\leq& n \left\Vert g( \cdot , \sigma, n)  \right\Vert_{\mathcal{L}_p} + \left\Vert q \right\Vert _{\mathcal{L}_p} \nonumber \\ &=& n\mathcal{L}_{p, n}(\sigma) + \left\Vert q \right\Vert_{\mathcal{L}_p} \nonumber \\ &\leq& n\mathcal{L}_{p, n}(\sigma) + N - n \nonumber,
	\end{eqnarray}
	which concludes the proof of the Lemma. \hspace*{\fill}\qed

\vspace{2mm}
Returning to the proof of Theorem 2, from (\ref{eq:24}) we see that $1 \leq m \leq N$ and $N-m \leq 2^d - 1$. Therefore applying Lemma \ref{lem:2} with $p=2$,
\begin{eqnarray}
N \mathcal{L}_{2, N}(\tilde{\sigma}) &\leq& m \mathcal{L}_{2, m}(\tilde{\sigma}) + N - m \nonumber \\ &\leq& N \mathcal{L}_{2, m}(\tilde{\sigma}) + 2^d - 1 \nonumber.
\end{eqnarray}
Consequently,
\begin{equation} \label{eq:31}
\mathcal{L}_{2, N}(\tilde{\sigma}) \leq \mathcal{L}_{2, m}(\tilde{\sigma}) + (2^d - 1)/N,
\end{equation}
and concluding from (\ref{eq:25}) and (\ref{eq:31}),
\begin{equation*}
	\mathcal{L}_{2, N}(\tilde{\sigma}) \leq C(d)F_n(\sigma) + (2^d-1)/N
\end{equation*} 
we obtain the required statement. \hspace*{\fill}\qed

\noindent
\paragraph{Proof of Theorem 4.} Let $\sigma_N$ be a finite sequence contained in $[0, 1)^d$, and let $\tilde{\sigma}_n$ denote a symmetric sequence consisting of $n = 2^dN$ terms, which is generated by $\sigma_N$.
Recall Theorem 1, 
\begin{equation} \label{eq:41}
\mathcal{L}_{2, n}(\tilde{\sigma}_n) \leq C(d)F_N(\sigma_N)
\end{equation}
and note the relation,
\begin{equation} \label{eq:42}
\frac{1}{C(d)} = 2^d\beta(d)
\end{equation}
with $\beta(d)$ defined as in (\ref{A3}). Therefore we can rewrite (\ref{eq:41}) as
\begin{equation} \label{eq:43}
F_N(\sigma_N) \geq \frac{1}{C(d)} \mathcal{L}_{2, n}(\tilde{\sigma}_n) = 2^{d}\beta(d) \mathcal{L}_{2, n}(\tilde{\sigma}_n),
\end{equation}
lower bounding the diaphony by the $\mathcal{L}_{2}-$discrepancy of the symmetrised sequence. At this time, we recall Roth's result regarding the lower bound for the $\mathcal{L}_2-$discrepancy. It states, 
\begin{equation*}
\mathcal{L}_{2, n}(\tilde{\sigma}) \geq \alpha(d-1)\frac{(\log n)^{\frac{d-1}{2}}}{n}
\end{equation*}
with constant $\alpha(d)$ defined as in (\ref{A2}). Thus, 
\begin{eqnarray} \label{eq:44}
\mathcal{L}_{2, n}(\tilde{\sigma}_n) &\geq& \alpha(d-1)n^{-1}(\log n)^{\frac{d-1}{2}} \nonumber \\ &=& \alpha(d-1)(2^dN)^{-1}\big(\log (2^dN)\big)^{\frac{d-1}{2}} \nonumber \\ &>& \alpha(d-1)2^{-d}N^{-1}(\log N)^{\frac{d-1}{2}}.
\end{eqnarray}
Putting together (\ref{eq:43}) and (\ref{eq:44}), we conclude that
\begin{eqnarray}
F_N(\sigma_N) &\geq& 2^d\beta(d)\mathcal{L}_{2, n}(\tilde{\sigma}_n) \nonumber \\ &>& 2^d\beta(d)[\alpha(d-1)2^{-d}N^{-1}(\log N)^{\frac{d-1}{2}}] \nonumber \\ &=& \alpha(d-1)\beta(d)N^{-1}(\log N)^{\frac{d-1}{2}} \nonumber
\end{eqnarray}
as required. \hspace*{\fill}\qed

\noindent
\paragraph{Proof of Theorem 5.}  Let $\sigma$ be an infinite sequence contained in $[0, 1)^d$. Let $\tilde{\sigma}$ be an infinite symmetric sequence contained in $[0, 1)^d$ which is generated by $\sigma$. Choose arbitrary constants $A(d), A_1(d)$ and $A_2(d)$ such that
\begin{equation*}
	0 < A(d) < A_1(d) < A_2(d) < \alpha(d)\beta(d).
\end{equation*}
Then clearly,
\begin{equation*}
	\frac{A_2(d)}{\beta(d)} < \alpha(d).
\end{equation*}
and from Theorem 3,
\begin{equation} \label{eq:46}
	N\mathcal{L}_{2, N}(\tilde{\sigma}) > \frac{A_2(d)}{\beta(d)}(\log N)^{\frac{d}{2}}
\end{equation} 
for infinitely many $N$. We choose sufficiently large enough $N$ which satisfies (\ref{eq:46}) and the conditions
\begin{equation} \label{eq:47}
	N \geq \frac{2^sA_1(d)}{A_1(d) - A(d)}
\end{equation}
\begin{center}
	\&
\end{center}
\begin{equation} \label{eq:48}
	A_2(d)(\log N)^{\frac{d}{2}} - (2^d - 1)\beta(d) \geq A_1(d)(\log N)^{\frac{d}{2}}.
\end{equation}
Observe from (\ref{eq:47}) that $N > 2^d$. Set $n = \lfloor N/2^d \rfloor$ and rearranging the statement of Theorem 2 with (\ref{eq:42}), we can write
\begin{equation} \label{eq:49}
	NF_n(\sigma) \geq 2^d\beta(d)N \mathcal{L}_{2, N}(\tilde{\sigma}) - 2^d(2^d - 1)\beta(d).
\end{equation}

\noindent
Using (\ref{eq:47}), 	
\begin{equation*}
	n = \lfloor N/2^d \rfloor > \frac{N}{2^d} - 1 \geq \bigg( \frac{A(d)}{A_1(d)} \bigg)  \bigg( \frac{N}{2^d} \bigg),
\end{equation*}
and conversely, 
\begin{equation*}
	n = \lfloor N / 2^d \rfloor \leq \frac{N}{2^d}.
\end{equation*}
Putting the last two inequalities together, we obtain
\begin{equation} \label{eq:50}
	n < N\leq\frac{2^dA_1(d)n}{A(d)}.
\end{equation}
From (\ref{eq:46}), (\ref{eq:48}) and (\ref{eq:49})
\begin{eqnarray}
	NF_n(\sigma) &\geq& 2^dA_2(d)(\log N)^{\frac{d}{2}} - 2^d(2^d - 1)\beta(d) \nonumber \\ &\geq& 2^dA_1(d)(\log N)^{\frac{d}{2}} \nonumber.
\end{eqnarray}
It follows from the last line and (\ref{eq:50}), that 
\begin{equation} \label{eq:51}
	F_n(\sigma) > A(d)n^{-1}(\log n)^{\frac{d}{2}}.
\end{equation}
The statement is now proved for all $n < N$. 

\vspace{2mm}
To show for all infinitely many $n \in \mathbb{N}$, we proceed via the following. Let $(N_k)_{k \in \mathbb{N}}$ be an infinite sequence of natural numbers satisfying the condition,
\begin{equation} \label{eq:52}
	N_{k+1} > N_k + 2^d \hspace{5mm} (k = 1, 2, ...)
\end{equation}
Let each member of the sequence $(N_k)_{k \in \mathbb{N}}$ also satisfy conditions (\ref{eq:46}) to (\ref{eq:48}) with $N = N_k, \hspace{1mm} (k = 1, 2, ...)$. For each $k \in \mathbb{N}$, set $n_k = \lfloor N_k / 2^d \rfloor$ and by (\ref{eq:52}) it follows that, 
\begin{equation*}
		n_1 < n_2 < n_3 < ...
\end{equation*}
With the same steps that we used to prove (\ref{eq:51}), we can conclude that 
\begin{equation*}
	F_{n_k}(\sigma) > A(d)n^{-1}_k(\log n_k)^{\frac{d}{2}}
\end{equation*}
for each $k \in \mathbb{N}$. Therefore, the estimate (\ref{eq:51}) is satisfied for infinitely many $k \in \mathbb{N}$. Given that the constant $A(d)$ is an arbitrary positive number less than $\alpha(d)\beta(d)$, we have the Theorem. \hspace*{\fill}\qed

\section{A Key Proof for the Dyadic Case}

In this final section we prove the main result, Theorem 8 and note that Theorems 9$-$11 can be shown along the same lines as the corresponding Theorems for the classical diaphony, just incorporating Theorem 8 instead of Theorem 1.

\noindent
\paragraph{Proof of Theorem 8.} Let $\sigma_n = (\mathbf{a}_i)_{i = 1}^n \subset [0, 1)^d$ and $\tilde{\sigma}_N = (\mathbf{b}_i)_{i=1}^N \subset [0, 1)^d$ be finite sequences consisting of $n$ and $N=2^dn$ terms respectively, such that $\sigma_n$ generates $\tilde{\sigma}_N$. Then from the definition of the sequence $\tilde{\sigma}_N$, we can rewrite the discrepancy function as follows.
\begin{equation} \label{eq:169}
	g \big( [\mathbf{0}, \boldsymbol{\gamma}), \tilde{\sigma}_N, N \big) = \frac{1}{2^dn} \sum_{i=1}^{n} \sum_{\boldsymbol{\theta} \in Z_d} \chi_{\boldsymbol{\gamma}} \bigg( \frac{1}{2}(\mathbf{1} - \boldsymbol{\theta}) + \boldsymbol{\theta} \mathbf{a}_i \bigg) - \lambda_d \big( [\mathbf{0}, \boldsymbol{\gamma}) \big),
\end{equation}
where $Z_d$ is defined as in Remark 1. The discrepancy function can then be written asymptotically equal to a Walsh series. That is,
	
\begin{equation*} 
	g \big( [\mathbf{0}, \boldsymbol{\gamma}), \tilde{\sigma}_N, N \big) \sim \sum_{\mathbf{k} \in \mathbb{N}_0^d} \widehat{g}(\mathbf{k}) \text{wal}_{\mathbf{k}}(\boldsymbol{\gamma})
\end{equation*}
where each of the Walsh coefficients $\widehat{g}(\mathbf{k})$ can be calculated by,	
\begin{equation} \label{walshcoeff}
	\widehat{g}(\mathbf{k}) = \int_{[0, 1)^d} g \big( [\mathbf{0}, \boldsymbol{\gamma}), \tilde{\sigma}_N, N \big) \text{wal}_{\mathbf{k}}(\boldsymbol{\gamma}) d\boldsymbol{\gamma}.
\end{equation}
 
 The following Lemma allows the computation of the integrals arising from the expression above for the Walsh coefficients.
 
\begin{lemma} \label{eq:lem12}
	Let $k \in \mathbb{N}$. For $x \in [0, 1)$, we define $a_n = a_n(x)$ and $b_n = b_n(x)$ by,
	\begin{equation}
	a_n := m \cdot 2^{-n} \leq x < (m+1) \cdot 2^{-n} =: b_n.
	\end{equation}
	for some integers $0 \leq m < 2^n$ and $n \geq 0$.
	\noindent
	Define and denote two integrals by,
	\begin{equation*}
	A(k, x) := \int_{x}^{1} \textnormal{wal}_k(\gamma) d\gamma  \hspace{5mm} \& \hspace{5mm}  B(k) := \int_{0}^{1} \gamma \hspace{0.5mm} \textnormal{wal}_{k}(\gamma) d\gamma.
	\end{equation*}
	Then,
	\begin{equation*}
	A(k, x) = \textnormal{wal}_{k}(x) \big( \psi_n - x \big) \hspace{8mm} \& \hspace{8mm} B(k) = 0 
	\end{equation*}
	where $\psi_n$ is defined as one of $a_n$ or $b_n$ depending on which is nearer to $x$. $($If $x$ is the midpoint of $(a_n, b_n)$, then set $\psi_n = b_n.)$
\end{lemma}

\noindent
\textit{Proof.}
We begin with $A(k, x)$. Section 3 of \cite{Fine} discusses integrals of the form,
\begin{equation*}
J_k(x) := \int_{0}^{x} \text{wal}_k(\gamma) d\gamma
\end{equation*}
and gives a concise and intuitive result on how one can compute integrals of this kind. Namely,
\begin{equation*}
J_k(x) = \text{wal}_k(x) \big( x - \psi_n \big)
\end{equation*}
where $\psi_n$ is defined as in the statement of the Lemma above. Due to an elementary fact in the study of Walsh functions,
\begin{equation*}
\int_{0}^{1} \text{wal}_k(\gamma) d \gamma = 0
\end{equation*}
for all $k \neq 0$ and we can therefore conclude as required,
\begin{eqnarray} \label{Aint}
A(k, x) = \int_{x}^{1} \textnormal{wal}_k(\gamma) d\gamma &=& \int_{0}^{1} \text{wal}_k(\gamma) d \gamma - \int_{0}^{x} \textnormal{wal}_k(\gamma) d\gamma \nonumber \\ &=& - J_k(x) \nonumber \\ &=& \text{wal}_k(x) \big( \psi_n - x \big).
\end{eqnarray}

\vspace{1mm}
\noindent
Moving on to the integral $B(k)$, we use a simple integration by parts and (\ref{Aint}) to show that
\begin{eqnarray}
B(k) = \int_{0}^{1} \gamma \hspace{0.5mm} \textnormal{wal}_{k}(\gamma) d\gamma \nonumber &=& - \int_{0}^{1} J_k(\gamma) d\gamma \nonumber \\ &=& \int_{0}^{1} \text{wal}_k(\gamma) \big( \psi_n - \gamma \big) d\gamma \nonumber \\ &=& \psi_n \int_{0}^{1} \text{wal}_k(\gamma) d\gamma - B(k) \nonumber.
\end{eqnarray}
Therefore, conclude that $B(k) = 0$. \hspace*{\fill} \qed

\vspace{1mm}
We now return to the main body of the proof of Theorem 8. Use the expression for the discrepancy function in (\ref{eq:169}) and set
\begin{equation*} 
	\boldsymbol{\zeta}_i(\boldsymbol{\theta}) \coloneqq \frac{1}{2} (\mathbf{1} - \boldsymbol{\theta}) + \boldsymbol{\theta} \mathbf{a}_i.
\end{equation*}
The Walsh coefficients from (\ref{walshcoeff}) become,
\begin{eqnarray} \label{walsh1}
	\widehat{g}(\mathbf{k}) &=& \frac{1}{2^dn} \sum_{i=1}^{n} \sum_{\boldsymbol{\theta} \in Z_d} \int_{[0, 1)^d} \chi_{\boldsymbol{\gamma}}(\boldsymbol{\zeta}_{i}) \text{wal}_{\mathbf{k}}(\boldsymbol{\gamma}) d\boldsymbol{\gamma} - \int_{[0, 1)^d} \lambda_d \big( [\mathbf{0}, \boldsymbol{\gamma}) \big) \text{wal}_{\mathbf{k}}(\boldsymbol{\gamma}) d\boldsymbol{\gamma} \nonumber \\ &=& \frac{1}{2^dn} \sum_{i=1}^{n} \sum_{\boldsymbol{\theta} \in Z_d} \int_{\boldsymbol{\zeta}_{i}}^{\mathbf{1}} \text{wal}_{\mathbf{k}}(\boldsymbol{\gamma}) d\boldsymbol{\gamma} - \int_{[0, 1)^d} \lambda_d \big( [\mathbf{0}, \boldsymbol{\gamma}) \big)  \text{wal}_{\mathbf{k}}(\boldsymbol{\gamma}) d\boldsymbol{\gamma} \nonumber \\ &=& \frac{1}{2^dn} \sum_{i=1}^{n} \sum_{\boldsymbol{\theta} \in Z_d} \prod_{j=1}^{d} \int_{\zeta_{ij}}^{1} \text{wal}_{k_j}(\gamma_j) d\gamma_j - \prod_{j=1}^{d} \int_{0}^{1} \gamma_j \hspace{0.5mm} \text{wal}_{k_j}(\gamma_j) d\gamma_j.
\end{eqnarray}
	
\noindent
From here we can first consider the coefficient $\widehat{g}(\mathbf{0})$, noting that $\text{wal}_{\mathbf{0}} = 1.$
\begin{eqnarray}
	\widehat{g}(\mathbf{0}) &=& \frac{1}{2^dn} \sum_{i=1}^{n} \sum_{\boldsymbol{\theta} \in Z_d} \prod_{j=1}^{d} \int_{\zeta_{ij}}^{1}  d\gamma_j - \prod_{j=1}^{d} \int_{0}^{1} \gamma_j d\gamma_j \nonumber \\ &=& \frac{1}{2^dn} \sum_{i=1}^{n} \prod_{j=1}^{d} \sum_{\theta = \pm 1} (1-\zeta_{ij}) - \prod_{j=1}^{d} \frac{1}{2} \nonumber \\ &=& \frac{1}{2^dn} \sum_{i=1}^{n} 1 - \frac{1}{2^d} = 0.
\end{eqnarray} 
	
\noindent
As mentioned in Section 2.1, Parseval's identity holds for the Walsh function system. Hence,
\begin{eqnarray} \label{eq:107}
	\mathcal{L}_{2, N}^2(\tilde{\sigma}_N) &=& \int_{[0, 1)^d} \Big| g \big( [\mathbf{0}, \boldsymbol{\gamma}), \tilde{\sigma}_N, N \big) \Big|^2 d\boldsymbol{\gamma} \nonumber \\ &=& \sum_{\mathbf{k} \in \mathbb{N}_0^d} \big| \widehat{g}(\mathbf{k}) \big|^2 = \sideset{}{'}\sum_{\mathbf{k} \in \mathbb{N}_0^d} \big| \widehat{g}(\mathbf{k}) \big|^2
\end{eqnarray}
where $\sideset{}{'}\sum$ denotes the summation without the zero index.

\vspace{1mm}
\noindent
For an arbitrary nonempty subset $A$ of the set $\{1, 2, ..., d\}$, denote by $M(A)$ the set consisting of points $\mathbf{k} = (k_1, ..., k_d) \in \mathbb{N}_0^d$ such that $k_j \neq 0 \hspace{1mm} (1 \leq j \leq d)$ if and only if $j \in A$. We define,
\begin{equation} \label{eq:108}
	\phi'(A) := \sum_{\mathbf{k} \in M(A)} \big| \widehat{g}(\mathbf{k}) \big| ^2.
\end{equation}
It is easy to see that (\ref{eq:107}) can be written in the form
\begin{eqnarray} \label{eq:109}
	\mathcal{L}_{2, N}^2(\tilde{\sigma}_N) &=& \sum_{p=1}^{d} \sum_{|A| = p} \sum_{\mathbf{k} \in M(A)} \big| \widehat{g}(\mathbf{k}) \big| ^2 \nonumber \\ &=& \sum_{p=1}^{d} \sum_{|A| = p} \phi'(A)
\end{eqnarray}
where the sum is over all subsets $A$ of $\{1, 2, ...., d\}$ such that $|A| = p$, for $1 \leq p \leq d$.

\vspace{1mm}
\noindent
Now fix $A$ to be some $p$ element subset of $\{1, 2, ..., d\}$. From Lemma \ref{eq:lem12} and taking $\mathbf{k} \in M(A)$, (\ref{walsh1}) simplifies to
	\begin{equation*}
	\widehat{g}(\mathbf{k}) = \frac{1}{2^dn} \sum_{i=1}^{n} \sum_{\boldsymbol{\theta} \in Z_d} \prod_{j \in A} \text{wal}_{k_j} (\zeta_{ij}) \big( \psi_{n_j} - \zeta_{ij} \big).
	\end{equation*}
	Then, noticing that $(\psi_{n_j} - \zeta_{ij}) \leq 2^{-n_j}$ when we write $k_j = 2^{n_j}+k_j'$ for integers $0 \leq k_j' < 2^{n_j}$ and $n_j \geq 0$ for each $j \in A$, we obtain
	\begin{equation} \label{eq:116}
	\widehat{g}(\mathbf{k}) \leq \frac{1}{2^dn} \sum_{i=1}^{n} \prod_{j \in A} \sum_{\theta_j = \pm 1} \text{wal}_{k_j} (\zeta_{ij}) \cdot 2^{-n_j}.
	\end{equation}

	\vspace{1mm}
	\noindent
	Considering the product in (\ref{eq:116}), first note the following elementary facts regarding the dyadic Walsh functions. The Walsh functions are periodic with period one, and therefore $\text{wal}_{k}(1-x) = \text{wal}_{k}(-x)$. Furthermore, for all $k \in \mathbb{N}_0$ and $x \in [0, 1)$, we have
	\begin{equation*}
	\textnormal{wal}_k(x) = \textnormal{wal}_k (- x).
	\end{equation*}
	Thus,
	\begin{eqnarray} \label{eq:117}
	\prod_{j \in A} \sum_{\theta_j = \pm 1} \text{wal}_{k_j} (\zeta_{ij}) \cdot 2^{-n_j} \nonumber &=& \prod_{j \in A} 2^{-n_j} \Big( \text{wal}_{k_j} (a_{ij}) + \text{wal}_{k_j} (1-a_{ij}) \Big) \\ &=& \prod_{j \in A} 2^{-n_j} \Big( \text{wal}_{k_j}(a_{ij}) + \text{wal}_{k_j}(- a_{ij}) \Big) \nonumber \\ &=& \prod_{j \in A} 2^{-n_j} \Big( 2 \hspace{0.25mm} \text{wal}_{k_j}(a_{ij}) \Big) \nonumber \\ &=& 2^p \hspace{0.5mm} \text{wal}_{\mathbf{k}}(\mathbf{a}_{i}) \prod_{j \in A} 2^{-n_j}.
	\end{eqnarray}
	
	\noindent
	 From (\ref{eq:116}) and (\ref{eq:117}),
	\begin{equation*}
	\widehat{g}(\mathbf{k}) \leq \frac{1}{2^{d-p}n} \prod_{j \in A} 2^{-n_j} \sum_{i=1}^{n} \text{wal}_{\mathbf{k}}(\mathbf{a}_{i})
	\end{equation*}
	thus we obtain the following estimate of the Walsh coefficients,
	\begin{equation} \label{estimateWalsh}
	\big| \widehat{g}(\mathbf{k}) \big| \leq \frac{1}{2^{d-p}} \prod_{j \in A} 2^{-n_j}  \Bigg| \frac{1}{n} \sum_{i=1}^{n} \text{wal}_{\mathbf{k}}(\mathbf{a}_{i}) \Bigg|.
	\end{equation}
	Squaring (\ref{estimateWalsh}),
	\begin{eqnarray} \label{eq:118}
	\big| \widehat{g}(\mathbf{k}) \big|^2 &\leq& \frac{1}{4^{d-p}} \prod_{j \in A} 2^{-2n_j} \Bigg| \frac{1}{n} \sum_{i=1}^{n} \text{wal}_{\mathbf{k}}(\mathbf{a}_{i}) \Bigg|^2 \nonumber \\ &=& \frac{r_2(\mathbf{k})}{4^{d-p}} \Bigg| \frac{1}{n} \sum_{i=1}^{n} \text{wal}_{\mathbf{k}}(\mathbf{a}_{i}) \Bigg|^2
	\end{eqnarray}
	where we have denoted $\prod_{j \in A} r_2(k_j)$ by $r_2(\mathbf{k})$, and the function $r_2$ is as in the definition of the dyadic diaphony. Using (\ref{eq:108}) and (\ref{eq:118}), we get
	\begin{eqnarray} \label{eq:119}
	\phi'(A) &\leq& \frac{1}{4^{d-p}} \sum_{\mathbf{k} \in M(A)} r_2(\mathbf{k}) \Bigg| \frac{1}{n} \sum_{i=1}^{n} \text{wal}_{\mathbf{k}}(\mathbf{a}_{i}) \Bigg|^2 \nonumber \\ &\leq& \sum_{\mathbf{k} \in M(A)} r_2(\mathbf{k}) \Bigg| \frac{1}{n} \sum_{i=1}^{n} \text{wal}_{\mathbf{k}}(\mathbf{a}_{i}) \Bigg|^2.
	\end{eqnarray}
	Now concluding from (\ref{eq:108}), (\ref{eq:109}) and (\ref{eq:119}),
	\begin{eqnarray}
	\mathcal{L}_{2, N}^2(\tilde{\sigma}_N) &\leq& \sum_{p=1}^{d} \sum_{|A| = p} \sum_{\mathbf{k} \in M(A)}  r_2(\mathbf{k}) \Bigg| \frac{1}{n} \sum_{i=1}^{n} \text{wal}_{\mathbf{k}}(\mathbf{a}_{i}) \Bigg|^2 \nonumber \\ &=& \sum_{\mathbf{k} \in \mathbb{N}_0^d \setminus \{\mathbf{0}\}}  r_2(\mathbf{k}) \Bigg| \frac{1}{n} \sum_{i=1}^{n} \text{wal}_{\mathbf{k}}(\mathbf{a}_{i}) \Bigg|^2 \nonumber \\ &=& \big( 3^d - 1 \big) \Bigg( \frac{1}{3^d - 1} \sum_{\mathbf{k} \in \mathbb{N}_0^d \setminus \{\mathbf{0}\}}  r_2(\mathbf{k}) \Bigg| \frac{1}{n} \sum_{i=1}^{n} \text{wal}_{\mathbf{k}}(\mathbf{a}_{i}) \Bigg|^2 \Bigg) \nonumber \\ &=& \delta^2(d) F^2_{2, n}(\sigma_n) \nonumber
	\end{eqnarray}
	with constant $\delta(d)$ is defined as in (\ref{C1}). \hspace*{\fill}\qed

\paragraph{Acknowledgements.}

\noindent
The author would like to express thanks in the first instance to Friedrich Pillichshammer for supplying the chapters of Proinov's monograph. Gratitude is also expressed to Florian Pausinger for the numerous, useful discussions in the development of this manuscript. A mention must be given to the anonymous referee for the valuable feedback which was used while improving the paper.

\pagebreak

\section*{An Appendix - Proof of Theorem A}

For the readers benefit we give a proof of Theorem A, the recent $\mathcal{L}_2-$discrepancy result used in Section 2.3 of our paper to improve the lower bound results of the one-dimensional constants.

\vspace{1mm}
The sketch of the proof, originally published in [10], was later found to contain a small inaccuracy and therefore after introducing some necessary preliminary material, a complete and rectified proof is provided.

\subsection*{Preliminaries}

We begin by noting that for the purposes of this appendix, an altered form of the discrepancy function is used as in [10]. Let $C_{\mathbf{z}} \coloneqq (z_1, 1] \times \dots \times (z_d, 1]$ for $\mathbf{z} = (z_1, \dots, z_d) \in [0, 1)^d$ and $\mathbf{x} \in [0, 1)^d$ be an arbitrary point. Then define,
\begin{equation*}
g(\mathbf{x}, \sigma_N, N) \coloneqq \sum_{\mathbf{z} \in \sigma_N} \chi_{C_{\mathbf{z}}} (\mathbf{x}) - N \lambda_d \big( [\mathbf{0}, \mathbf{x}) \big)
\end{equation*}
for a finite $N$ term sequence, $\sigma_N \subset [0, 1)^d$. Notice that the summation in the definition above is simply the number of terms of the sequence $\sigma_N$ that are contained in the subinterval $[\mathbf{0}, \mathbf{x})$.

\paragraph{Haar$-$coefficients of the Discrepancy Function.} A dyadic interval of length $2^{-j} \hspace{1mm} (j \in \mathbb{N}_0)$ in $[0, 1)$, is an interval of the form 
\begin{equation*}
I \coloneqq \bigg[ \frac{m}{2^j}, \frac{m+1}{2^j} \bigg)
\end{equation*} 
for $m = 0, 1, \dots , 2^j - 1$. The Haar function $h_I = h_{j, m}$ with support $I$ is the function on $[0, 1)$ which is $1$ on the left half of $I$ , $-1$ on the right half of $I$, and $0$ outside of $I$. The $\mathcal{L}_{\infty}-$normalised Haar system consists of all Haar functions $h_{j, m}$ (for $j \in \mathbb{N}_0$ and $m = 0, 1, \dots , 2^j - 1$) together with the indicator function of $[0, 1)$, $h_{-1, 0}$. After normalisation in $\mathcal{L}_2\big([0, 1) \big)$, we obtain the orthonormal Haar basis of $\mathcal{L}_2 \big([0, 1) \big)$.

\vspace{2mm}
Let $\mathbb{N}_{-1} = \{-1, 0, 1, 2, \dots\}$, and define $\mathbb{D}_j = \{0, 1, \dots , 2^j - 1\}$ for $j \in \mathbb{N}_{0}$ and $\mathbb{D}_{-1} = \{0\}$ for $j = -1$. In higher dimensions, i.e. $d \geq 2$, the Haar function $h_{\mathbf{j}, \mathbf{m}}$ is given as the tensor product $h_{\mathbf{j}, \mathbf{m}}(\mathbf{x}) = h_{j_1, m_1}(x_1) \dots h_{j_d, m_d}(x_d)$ for $\mathbf{x} = (x_1, \dots , x_d) \in [0, 1)^d, \mathbf{j} = (j_1, \dots , j_d) \in \mathbb{N}_{-1}^d$ and $\mathbf{m} = (m_1, \dots , m_d) \in \mathbb{D}_{\mathbf{j}} \coloneqq \mathbb{D}_{j_1} \times \dots \times \mathbb{D}_{j_d}$. We will call the subintervals of $[0, 1)^d$, $I_{\mathbf{j}, \mathbf{m}} = I_{j_1, m_1} \times \dots \times I_{j_d, m_d}$ dyadic boxes. Note that all dyadic boxes with fixed $\mathbf{j}$ are congruent, hence we call $\mathbf{j}$ the shape of the box $I_{\mathbf{j}, \mathbf{m}}$. Lastly for $\mathbf{j} \in \mathbb{N}_{-1}^d$, let $|\mathbf{j}| = \max(0, j_1) + \max(0, j_2) + \dots + \max(0, j_d)$. The $\mathcal{L}_{\infty}-$normalised tensor Haar system consists of all Haar functions $h_{\mathbf{j}, \mathbf{m}}$ with $\mathbf{j} \in \mathbb{N}_{-1}^d$ and $\mathbf{m} \in \mathbb{D}_{\mathbf{j}}$. After normalisation in $\mathcal{L}_2 \big([0, 1)^d \big)$, we obtain the orthonormal Haar basis of $\mathcal{L}_2 \big([0, 1)^d\big)$.

\vspace{2mm}
Parseval's equality shows that the $\mathcal{L}_2-$norm of a function $f \in \mathcal{L}_2 \big( [0, 1)^d \big)$, denoted $\norm{f}_{\mathcal{L}_2}$, can be computed as 
\begin{equation*} \label{Haarparseval}
\norm{f}^2_{\mathcal{L}_2} = \sum_{\mathbf{j} \in \mathbb{N}_{-1}^d} 2^{|\mathbf{j}|} \sum_{\mathbf{m} \in \mathbb{D}_{\mathbf{j}}} \big| \mu_{\mathbf{j}, \mathbf{m}} \big|^2, 
\end{equation*}
where
\begin{equation*} \label{Haarcoeff}
\mu_{\mathbf{j}, \mathbf{m}} = \mu_{\mathbf{j}, \mathbf{m}}(f) = \int_{[0, 1)^d} f(\mathbf{x}) \cdot h_{\mathbf{j}, \mathbf{m}}(\mathbf{x}) d\mathbf{x}
\end{equation*}
are the Haar-coefficients of $f$.

\subsection*{Useful Lemmas}

We will require the following Lemmas.

\begin{statement}{Lemma A}
	Let $f(\mathbf{x}) = x_1 \dots x_d$ for $\mathbf{x} = (x_1, \dots, x_d) \in [0, 1)^d$. Let $\mathbf{j} \in \mathbb{N}^d_{0}, \mathbf{m} \in \mathbb{D}_{\mathbf{j}}$, and let $\mu_{\mathbf{j}, \mathbf{m}}$ be the Haar-coefficient of $f$. Then, 
	\begin{equation*}
	\mu_{\mathbf{j}, \mathbf{m}} = 2^{-2|\mathbf{j}|-2d}.
	\end{equation*}
\end{statement}

	\noindent
	\textit{Proof.}
	For $\mathbf{x} \in [0, 1)^d$, $\mathbf{j} \in \mathbb{N}_0^d$ and $\mathbf{m} \in \mathbb{D}_j$, let 
	\begin{eqnarray}
	\mu_{\mathbf{j}, \mathbf{m}} &=& \int_{[0, 1)^d} \lambda_d \big( [\mathbf{0}, \mathbf{x}) \big) \cdot h_{\mathbf{j}, \mathbf{m}}(\mathbf{x}) d\mathbf{x} \nonumber \\ &=& \int_{[0, 1)^d} x_1 \dots x_d \cdot h_{\mathbf{j}, \mathbf{m}}(\mathbf{x}) d\mathbf{x} \nonumber \\ &=& \prod_{i=1}^d \int_{0}^{1} x_i \cdot h_{j_i, m_i} (x_i) dx_i. \nonumber
	\end{eqnarray}
	
	\noindent
	Note that for $j_i \in \mathbb{N}_0$,
	\begin{equation*}
	\int_{0}^{1} x_i \cdot h_{j_i, m_i} (x_i) dx_i = 2^{-2j_i-2}.
	\end{equation*}
	
	\noindent
	Therefore, we can conclude easily that
	\begin{eqnarray}
	\mu_{\mathbf{j}, \mathbf{m}} &=& \prod_{i=1}^{d} 2^{-2j_i -2} \nonumber \\ &=& 2^{-2j_1 - 2 - \dots - 2j_d -2} \nonumber \\ &=& 2^{-2|\mathbf{j}|-2d} \nonumber
	\end{eqnarray}
	as required. \hspace*{\fill}\qed

\begin{statement}{Lemma B}
	Fix $\mathbf{z} = (z_1, \dots, z_d) \in [0, 1)^d$, and let $f(\mathbf{x}) = \chi_{C_{\mathbf{z}}}(\mathbf{x})$ be the characteristic function for the subinterval $C_{\mathbf{z}} = (z_1, 1] \times \dots \times (z_d, 1]$ with $\mathbf{x} = (x_1, \dots, x_d) \in [0, 1)^d$ arbitrary. Let $\mathbf{j} \in \mathbb{N}^d_{0}, \mathbf{m} \in \mathbb{D}_{\mathbf{j}}$ and $\mu_{\mathbf{j}, \mathbf{m}}$ be the Haar-coefficient of $f$. Then
	\begin{equation*}
	\mu_{\mathbf{j}, \mathbf{m}} = 0,
	\end{equation*}
	whenever $\mathbf{z}$ is not contained in the interior of the dyadic box $I_{\mathbf{j}, \mathbf{m}}$ supporting $h_{\mathbf{j}, \mathbf{m}}$.
\end{statement}

	\noindent
	\textit{Proof.}
	Take $\mathbf{z} \in [0, 1)^d$, such that $\mathbf{z} \notin I_{\mathbf{j}, \mathbf{m}}$. Note in the first instance, that
	\begin{eqnarray}
	\mu_{\mathbf{j}, \mathbf{m}} &=& \int_{[0, 1)^d} \chi_{C_{\mathbf{z}}}(\mathbf{x}) \cdot h_{\mathbf{j}, \mathbf{m}} (\mathbf{x}) \hspace{0.5mm} d\mathbf{x} \nonumber \\ &=& \int_{\mathbf{z}}^{\mathbf{1}} h_{\mathbf{j}, \mathbf{m}} (\mathbf{x}) \hspace{0.5mm} d\mathbf{x} \nonumber \\ &=& \prod_{i=1}^{d} \int_{z_i}^{1} h_{j_i, m_i} (x_i) dx_i. \nonumber
	\end{eqnarray}
	
	\noindent
	So $\mathbf{z} \notin I_{\mathbf{j}, \mathbf{m}}$ implies that $z_i \notin I_{j_i, m_i} := \Big[ \frac{m_i}{2^{j_i}}, \frac{m_i+1}{2^{j_i}} \Big)$ for some $i \in \{1, 2, \dots, d\}$ and $m_i \in \mathbb{D}_{j_i}$. Thus, there are two cases to consider. Either,
	\begin{enumerate}
		\item $z_i \geq \frac{m_i+1}{2^{j_i}}$, in which case $h_{j_i, m_i} = 0$ for all $\frac{m_i+1}{2^{j_i}} \geq x \geq 1$. This implies, 
		\begin{equation*}
		\int_{z_i}^{1} h_{j_i, m_i}(x_i) dx_i = 0.
		\end{equation*}
		Or, 
		\item $z_i < \frac{m_i}{2^{j_i}}$, in which case
		\begin{equation*}
		\int_{z_i}^{1} h_{j_i, m_i}(x_i) dx_i = 0
		\end{equation*}
		since $h_{j_i, m_i}$ equals $0$ in $I_{j_i, m_i}$ and $\int_{I_{j_i, m_i}} h_{j_i, m_i}(x_i) dx_i = 0$.
	\end{enumerate}
	
	\noindent
	Therefore, we can conclude for $\mathbf{z} \notin I_{\mathbf{j}, \mathbf{m}}$
	\begin{equation*}
	\mu_{\mathbf{j}, \mathbf{m}} = \int_{[0, 1)^d} \chi_{C_{\mathbf{z}}} (\mathbf{x}) \cdot h_{\mathbf{j}, \mathbf{m}}(\mathbf{x}) d\mathbf{x} = 0
	\end{equation*}
	as required. \hspace*{\fill}\qed

\subsection*{Main Statement \& Proof}

\begin{statement}{Theorem A}
	For a finite sequence $\sigma_N$ contained in $[0, 1)^d$, the inequality 
	\begin{equation*}
	\mathcal{L}_{2, N}(\sigma_N) \geq \gamma(d-1) \frac{(\log N)^{\frac{d-1}{2}}}{N}
	\end{equation*}
	holds with 
	\begin{equation*}
	\gamma(d) \coloneqq \frac{1}{\sqrt{21} \cdot 2^{2d+1} \sqrt{d!} (\log 2)^{\frac{d}{2}}}.
	\end{equation*}
\end{statement}

	\noindent
	\textit{Proof.}
	Let $\sigma_N$ be a finite sequence contained in $[0, 1)^d$, and take arbitrary $\mathbf{x} \in [0, 1)^d$. Let $\mathbf{j} \in \mathbb{N}^d_{0}, \mathbf{m} \in \mathbb{D}_{\mathbf{j}}$ be such that no point of $\sigma_N$ lies in the interior of the dyadic box $I_{\mathbf{j}, \mathbf{m}}$ supporting $h_{\mathbf{j}, \mathbf{m}}$. Let $\mu_{\mathbf{j}, \mathbf{m}}$ denote the Haar-coefficient of the discrepancy function, which for the purposes of this proof as already mentioned, we define as
	\begin{equation*}
	g(\mathbf{x}, \sigma_N, N) \coloneqq \sum_{\mathbf{z} \in \sigma_N} \chi_{C_{\mathbf{z}}} (\mathbf{x}) - N \lambda_d \big( [\mathbf{0}, \mathbf{x}) \big).
	\end{equation*}
	Now, the Lemmas above imply that
	\begin{equation*}
	\mu_{\mathbf{j}, \mathbf{m}} = -N 2^{-2|\mathbf{j}|-2d},
	\end{equation*}
	where $|\mathbf{j}| = j_1+j_2 + \dots + j_d.$ Next note that for a fixed $\mathbf{j} \in \mathbb{N}^d_{0}$, the cardinality of $\mathbb{D}_{\mathbf{j}}$ is $2^{|\mathbf{j}|}$, and the interiors of the dyadic boxes $I_{\mathbf{j}, \mathbf{m}}$ supporting $h_{\mathbf{j}, \mathbf{m}}$ are mutally disjoint. This implies that there are at least $2^{|\mathbf{j}|} - N$ such $\mathbf{m} \in \mathbb{D}_{\mathbf{j}}$ for which no point in the $N$ term sequence $\sigma_N$ lies in the interior of the dyadic box $I_{\mathbf{j}, \mathbf{m}}$ supporting $h_{\mathbf{j}, \mathbf{m}}$.
	
	\vspace{1mm}
	\noindent
	Set $M \coloneqq \lceil \log_2 N \rceil$. Then from Parseval's equality,
	\begin{eqnarray} \label{eqn5}
	\mathcal{L}^2_{2, N} (\sigma_N) &\geq& N^2 \sum_{|{\mathbf{j}}| \geq M} 2^{|{\mathbf{j}}|} \big( 2^{|{\mathbf{j}}|} - N \big) 2^{-4|\mathbf{j}| - 4d} \nonumber \\ &=& 	2^{-4d} N^2 \sum_{|{\mathbf{j}}| \geq M} 4^{-|{\mathbf{j}}|} - 2^{-4d} N^3 \sum_{|{\mathbf{j}}| \geq M} 8^{-|{\mathbf{j}}|} \nonumber
	\end{eqnarray}
	where the inequality is due to summing only those $\mu_{\mathbf{j}, \mathbf{m}}$ with $\mathbf{j} \in \mathbb{N}_0^d$ and $\mathbf{m} \in \mathbb{D}_j$ as chosen above. Considering the summations, the coefficient of $M^{d-1}$ in
	\begin{equation*}
	\sum_{|\mathbf{j}| \geq M} q^{-|\mathbf{j}|}
	\end{equation*}
	is computed as,
	\begin{equation*}
	\frac{q^{-M+1}}{(q-1)(d-1)!}
	\end{equation*}
	for $q>1$. This implies,
	\begin{eqnarray}
	\mathcal{L}^2_{2, N} (\sigma_N) &\geq& 2^{-4d} \big( N2^{-M} \big)^2 \frac{4}{3} \frac{M^{d-1}}{(d-1)!} - 2^{-4d} \big( N2^{-M} \big)^3 \frac{8}{7} \frac{M^{d-1}}{(d-1)!}. \nonumber
	\end{eqnarray}
	
	\noindent
	Now let $t = M - \log_2 N$, so that $0 \leq t < 1$ and $N2^{-M} = 2^{-t}$. Then we get,
	\begin{equation*}
	\mathcal{L}^2_{2, N} (\sigma_N) \geq \omega (\log_2 N)^{d-1}
	\end{equation*}
	if
	\begin{equation*}
	2^{-4d} \big( 2^{-2t} \big) \frac{4}{3} \frac{M^{d-1}}{(d-1)!} - 2^{-4d} \big( 2^{-3t} \big) \frac{8}{7} \frac{M^{d-1}}{(d-1)!} \geq \omega (M-t)^{d-1} 
	\end{equation*}
	which is satisifed if,
	\begin{equation*}
	\omega \leq \frac{1}{2^{4d}(d-1)!} \bigg( \frac{4}{3} \big(2^{-2t} \big) - \frac{8}{7}  \big( 2^{-3t} \big) \bigg)
	\end{equation*}
	for all $0 \leq t < 1$. Or equivalently,
	\begin{equation*}
	\omega \leq \frac{1}{2^{4d}(d-1)!} \bigg( \frac{4}{3} y^2 - \frac{8}{7} y^3 \bigg)
	\end{equation*}
	for all $\frac{1}{2} < y \leq 1$.
	
	\vspace{2mm}
	To finish, we require the minimal value of the expression above. This occurs when $y=1$, and equals $\frac{4}{21}$. It follows that we have arrived at the desired constant. To obtain the exact order in the statement, recall that we must divide by $N$ to rectify the use of the altered discrepancy function throughout the proof. \hspace*{\fill}\qed

\end{document}